\documentclass[12pt]{article}
\usepackage{amssymb}
\usepackage{mathrsfs}
\usepackage{amsfonts}
\usepackage{graphicx}
\usepackage{colortbl,dcolumn}
\usepackage{latexsym}
\usepackage{amsbsy,amsthm}
\usepackage{amsmath}
\usepackage{bbm}
\usepackage{psfrag}
\usepackage{float} 
\usepackage{exscale}
\usepackage{relsize}
\usepackage{booktabs}
\usepackage{verbatim}
\usepackage{enumerate}
\usepackage{cite}
\usepackage{amsthm}
\usepackage[top=0.8in, bottom=0.8in, left=0.8in, right=0.8in]{geometry}
\usepackage{hyperref}

\numberwithin{equation}{section}

\newtheorem{assumption}{Assumption}[section]

\newtheorem{lemma}[assumption]{Lemma}
\newtheorem{theorem}[assumption]{Theorem}   %共享编号

\newtheorem{corollary}[assumption]{Corollary}

\hypersetup{colorlinks=true,                         %给引用加超链接
linkcolor=blue,
anchorcolor=blue,
citecolor=blue}

\allowdisplaybreaks[4] %公式跨页
\begin{document}
\title
{
    Higher order numerical methods for SDEs without globally monotone coefficients\footnote
    {
    This work is supported by Natural Science Foundation of China (12471394, 12071488, 12371417), 
    %Natural Science Foundation of Hunan Province (2020JJ2040), 
    Postgraduate Scientific Research Innovation Project of Hunan Province (No: CX20230349) and Fundamental Research Funds for the Central Universities of Central South University (No: 2021zzts0480).
    }
}

\author
{
    Lei Dai, Xiaojie Wang\footnote
    {
    Corresponding author: x.j.wang7@csu.edu.cn; x.j.wang7@gmail.com.
    }
    \\
    \footnotesize  School of Mathematics and Statistics, HNP-LAMA, Central South University, Changsha, China
    \\
}

%\date{}
\maketitle
\begin{abstract}
{
    In the present work, 
%we design and analyze higher order strong approximations of SDEs with non-globally monotone coefficients. 
%This work is mainly inspired by [Hutzenthaler, Jentzen, Wang, Math. Comp., 2018] and [Hutzenthaler, Kisker, arXiv, 2022], and it significantly improves convergence rate of numerical methods in comparison to these works.
%More specifically, 
we delve into further study of numerical approximations of SDEs with non-globally monotone coefficients. We design and analyze a new family of stopped increment-tamed time discretization schemes of Euler, Milstein and order $1.5$ type for such SDEs.  
By formulating a novel unified framework, the proposed methods are shown to possess the exponential integrability properties,
which are crucial to recovering convergence rates in the non-global monotone setting. 
Armed with such exponential integrability properties and by the arguments of perturbation estimates, we successfully identify the optimal strong convergence rates of the aforementioned methods in the non-global monotone setting.  Numerical experiments are finally presented to corroborate the theoretical results.
}
%Moreover, this numerical method succeeds an essential feature, 
%ultimate boundedness property, of the original model.  
%we improved a series of results in the aforementioned reference and answered a question ([Remark 3.1]) raised by that paper, i.e., can stopped increment-tamed Euler-Marayuma method hold order one strong convergence rate for additive noise models? 

\end{abstract}

\section{Introduction}

Over last decades, there has been considerable interest in developing and analyzing time discretization schemes for stochastic differential equations (SDEs) that are extensively used in various fields of science. This interest is driven by the fact that 
%SDEs play an essential role in various fields of science and engineering, but 
the closed-form solutions of nonlinear SDEs are rarely available and one turns to reliable numerical approximations. Under globally Lipschitz conditions, the classical monographs \cite{1992Numerical,Milstein2004} established fundamental theory to analyze various numerical methods for SDEs, such as the classic Euler--Maruyama (EM) method and Milstein method. Nevertheless, a variety of SDEs in applications do not obey such a traditional but restrictive condition. One might wonder whether the commonly used schemes in the globally Lipschitz setting still perform well for SDEs with non-globally Lipschitz coefficients.
In 2011, the authors of \cite{hutzenthaler2011strong} gave a negative answer and
showed that the usual EM method is divergent in the case of SDEs with super-linearly growing coefficients 
that violate the globally Lipschitz condition. 

To discretize such SDEs with super-linear coefficients, a natural candidate is the implicit methods such as the backward Euler/Milstein methods, which have been extensively investigated, see, e.g., \cite{higham2002strong,beyn2016stochastic,wang2020mean,wang2023mean} and references therein. The other choice is to
design some explicit schemes based on modifications of traditional explicit EM/Milstein methods \cite{hutzenthaler2012strong,sabanis2016euler,Li2018explicit,Fang2020,kelly2017adaptive,Kelly2023}.
A typical example is the "tamed Euler" (TE) method, 
which is explicit and originally proposed by \cite{hutzenthaler2012strong}. 
%which has been shown to possess strong convergence rate of order $\tfrac{1}{2}$, known as  to approximate SDEs with the one-sided Lipschitz condition 
Since then, there has been a flourishing literature on numerical approximations of SDEs with non-globally Lipschitz coefficients, see, e.g., \cite{higham2012convergence,Tretyakov2013fundamental,wang2013tamed,neuenkirch2014first,mao2015truncated,beyn2016stochastic,sabanis2016euler,hutzenthaler2018exponential,Li2018explicit,Fang2020,brehier2023approximation,Kelly2023}. 
In order to obtain the desired convergence rates of numerical schemes, 
the aforementioned publications essentially relied on the popular global monotonicity condition: for some $q >1$,
    \begin{equation}\label{eq:oneside_coef_condi}
        \begin{aligned}
            \langle x - y,f(x) - f(y)\rangle  + \tfrac{q}{2}\|g(x) - g(y)\|^2 \le {K}|x - y{|^2}
            \quad
            \forall x,y \in \mathbb{R}^d,
%            \\
%            &\langle x,f(x)\rangle  + \tfrac{1}%{2}\|g(x)\|^2 \le {K_2}{(1 + |x|)^2},
        \end{aligned}
    \end{equation}
    where $f(\cdot)$ and $g(\cdot)$ are the drift and diffusion coefficients of SDEs, respectively and 
    $K > 0$ is a uniform constant independent of $x,y$.
Nevertheless, 
%the one-sided Lipschitz condition 
%    \begin{equation}\label{eq:oneside_coef_condi}
 %       \begin{aligned}
 %           \langle x - y,f(x) - f(y)\rangle   \leq {K}|x - y{|^2}, \|g(x)-g(y)\| \leq {K}|x - y|,\ K>0,
 %       \end{aligned}
 %   \end{equation}
%where $f(\cdot)$ and $g(\cdot)$ are the drift and diffusion coefficients of SDEs, respectively. 
%%%%%%%%%%
%This result significantly broadens the classical framework of approximating SDEs with globally Lipschitz coefficients, as presented in \cite{1992Numerical,Milstein2004}. 
%%%%%%%%%%
%or 
such a global monotonicity condition is still too restrictive and many practical SDEs cannot satisfy it, such as the stochastic Lorenz equation, the stochastic Langevin dynamics, the stochastic Lotka-Volterra model and so on (see, e.g., \cite{hutzenthaler2015numerical,mao2007stochastic,cox2024local}). 
In the article \cite{hutzenthaler2018exponential}, the authors proposed  an explicit "stopped increment-tamed Euler--Maruyama" (SITEM) method
for SDEs with non-globally monotone coefficients and established exponential integrability properties and strong convergence of the numerical approximations, but without a convergence rate revealed.
A recent breakthrough by \cite{Hutzenthaler2020} developed a framework of perturbation theory
that is crucial to recovering strong convergence rates of numerical schemes in a non-globally monotone setting.
More precisely, by combining  the perturbation theory with exponential integrability properties of the numerical approximations established in \cite{hutzenthaler2018exponential}, 
one could acquire a strong convergence rate of at most order $\tfrac{1}{2}$ for the SITEM method in certain non-globally monotone setting. 

An interesting question was raised by  \cite[Remark 3.1]{Hutzenthaler2020}, as to whether it is possible to achieve a higher order convergence rate (greater than $\tfrac12$) in the absence of global monotone setting. 
Unfortunately, following \cite[Theorem 1.2]{Hutzenthaler2020}, the convergence rates of any schemes would not exceed order $\tfrac12$, 
which is exactly the order of the H\"{o}lder regularity in the $L^p$ sense of the approximation process.
Very recently, the authors of \cite{dai2023order} partially answered the interesting question, %by presenting some novel perturbation estimates, which successfully 
by recovering the order-one convergence rate of the SITEM method in \cite{hutzenthaler2018exponential} for some particular SDEs such as a class of second order SDEs and additive noise driven SDEs. This was achieved by developing some new perturbation estimates and carry out more careful error estimates.

Still, it remains unclear whether order one (or higher order) strong convergence can be derived for approximations of more general SDEs with multiplicative noises in the non-globally monotone setting.
%
%Naturally, the next important topic is how to develop a higher order numerical method for the non-global monotonic SDEs, which is also the main issue of this paper. 
%
In this article, inspired by \cite{hutzenthaler2018exponential,hutzenthaler2022stopped,dai2023order}, we attempt to resolve the tough problem and introduce a class of new higher order numerical methods \eqref{eq:new_stop_tamed_type_methods} for general SDEs with multiplicative noises and non-global monotonic coefficients. 
By choosing different coefficients in \eqref{eq:new_stop_tamed_type_methods}, we obtain stopped increment-tamed time discretization schemes of Euler, Milstein and order $1.5$ type.
%such as the stopped increment-tamed Euler–Maruyama method , the stopped increment-tamed Milstein method and the stopped increment-tamed order $1.5$ method
%More precisely, we design a new family of stopped increment-tamed methods with exponential integrability properties. 
%
%Within a unified framework, the newly proposed methods are shown to possess the exponential integrability properties (see Corollary \ref{cor:numerical-exponential-moment}), which are crucial to recovering convergence rates in the non-global monotone setting. 
%
By formulating an exponential integrability theorem (see Theorem \ref{thm:ex_inte_property_of_stop_tamed_method}), one can easily derive the exponential integrability properties of these methods within a unified framework. 
Armed with such exponential integrability properties and by the arguments of perturbation estimates, %presented in \cite{dai2023order},  
we successfully reveal the pathwise uniformly order-one strong convergence rate for the stopped increment-tamed Milstein method applied to general SDEs (Theorem \ref{thm:conver_rate_stopped_tamed_milstein_method}). Following the same idea but with more delicate perturbation estimates, we derive the optimal strong convergence rate of the order $1.5$ method (Theorem \ref{thm:conver_rate_stopped_tamed_1.5_order_method}). To the best of our knowledge, this is the first work to recover order $1$ and $1.5$ of numerical method for general SDEs with non-global monotone coefficients.

%,  the optimal strong convergence rates of the Milstein and order $1.5$ type schemes in the non-global monotone setting (Theorems \ref{thm:conver_rate_stopped_tamed_milstein_method}, \ref{thm:conver_rate_stopped_tamed_1.5_order_method}).
%
%By selecting distinct tamed terms, we are able to obtain different numerical schemes such as the stopped increment-tamed Euler–Maruyama method (new), the stopped increment-tamed Milstein method and the stopped increment-tamed order $1.5$ method. 
%
%

It is worth mentioning that our contribution is highly non-trivial, compared to existing works.
Essential difficulties are twofold. First, developing higher order numerical schemes with exponential integrability properties is not an easy task. Secondly, obtaining the expected convergence rates greater than $\tfrac12$ is also much involved.
%
%The first challenge in devising higher order schemes lies in the preservation of exponential integrability properties. So far, we are aware of only a few references dealing with the exponential integrability properties of numerical approximations for nonlinear SDEs.
%By using truncation and taming strategy, the authors of \cite{hutzenthaler2018exponential} and \cite{hutzenthaler2019strong} constructed a class of numerical methods possessing exponential integrability properties for, respectively, the certain SDEs and the stochastic Burgers equation with non-global monotone coefficients. 
%Furthermore, the authors of \cite{hutzenthaler2022stopped} proposed a Brownian-increment tamed Euler method that can inherit exponential integrability properties. This method tames the increment of the Brownian motion and employs the same truncated technique as \cite{hutzenthaler2018exponential}. Another interesting perspective comes from the preservation of Hamiltonian (energy).  In \cite{cui2022density}, the authors developed an energy-preserving averaged vector field (AVF) scheme, which preserves  exponential integrability properties, to approximate the stochastic Langevin equation.
Indeed, we are only aware of a few works dealing with the exponential integrability properties of numerical approximations of nonlinear SDEs \cite{hutzenthaler2018exponential,hutzenthaler2019strong,hutzenthaler2022stopped,cui2022density}, where just Euler-type schemes were examined. However, the idea of obtaining exponential integrability properties there cannot be directly generalised to higher order methods and we introduce new ideas to overcome the difficulties.
%It is worth noting that all the previous results are Euler-type numerical methods and have a lower convergence rate. Moreover, the exponential integrability properties from these works cannot be generalised to higher order methods directly. For example, if the Milstein method is employed as the tamed term in the SITEM method, as opposed to the Euler–Maruyama method, then the conditions, such as equation $(56)$ in \cite[Lemma 2.8]{hutzenthaler2018exponential}, will not be satisfied.
%

%In the present paper, 
%inspired by \cite{hutzenthaler2018exponential} and \cite{hutzenthaler2022stopped}, 

The remainder of this paper is organised as follows. The next section introduces some notations and inequalities that may be used subsequently. In Section 3, we propose a new family of new methods and present their properties. Then the exponential integrability properties of these methods are elaborated in Section 4. Equipped with the exponential integrability properties, we are able to derive optimal strong convergence rates of these methods in Section 5. 
Finally, applications and numerical experiments are provided in Section 6.

\section{Preliminaries and notations}

Let $(\Omega,\mathcal{F},$\ $\{\mathcal{F}_t\}_{t\geq 0},\mathbb{P})$ be a complete probability space with a usual filtration $\{\mathcal{F}_t\}_{t\geq 0}$ and $\{W_t\}_{t\geq0}$ be an $m$-dimensional ($m \in \mathbb{Z}^{+}$) standard Brownian motion.
Consider the following It\^{o} SDEs:
    \begin{equation} \label{eq:typical_sde}
      \left\{ 
        \begin{array}{l}
        \displaystyle{
        X_t-X_0
        =
            \int_0^t f(X_s){\rm d}s 
                    +
                    \int_0^t g(X_s){\rm d}W_s,
                    \ 0 < t \leq T< +\infty},
        \\
        X_0 = \xi_X,
        \end{array} \right.
    \end{equation}
where $f: \mathbb{R}^d \rightarrow\mathbb{R}^d, d\in \mathbb{Z}^{+}$ stands for the drift coefficient, $g:\mathbb{R}^d \rightarrow\mathbb{R}^{d\times m}$ the diffusion coefficient and $\xi_X$ the initial data.
For $x\in \mathbb{R}^d$, the $x^{(i)},i=1,...,d$ denotes the $i$-th component of $x$ and $|x|$ denotes the Euclidean norm induced by the vector inner product $\langle\cdot,\cdot\rangle$.
%%%
For a matrix  $A \in \mathbb{R}^{d\times m}$, $A^{(i)},i=1,...,m$ denotes the $i$-th column of $A$ and $A^{(ij)},i=1,...,d,j=1,...,m$ represents the element at $i$-th row and $j$-th column of $A$.
Let $A^*$ be the transpose of $A$ and $\|A\|:= \sqrt{\operatorname{trace}(A^{*}A)}$ be the Hilbert-Schmidt norm induced by the Hilbert-Schmidt inner product $\langle\cdot,\cdot\rangle_{HS}$.
For $
    \Lambda : \mathbb{R}^d
    \rightarrow
    \mathbb{R}^{z_1 \times z_2}, z_1,z_2\in \mathbb{Z}^+ 
    $ 
and $y\in \mathbb{R}^d$,
we denote  
\begin{equation}
    \begin{aligned}
        \mathscr{L}^j_g \Lambda(y)
        &: = 
            \sum\limits_{i=1}^d 
                g^{(ij)} (y)
                \tfrac{\partial \Lambda}
                      {\partial x^{(i)}}(y),\ j=1,...,m,
    \end{aligned}
\end{equation}
and
\begin{equation}
    \begin{aligned}
        \mathcal{A}_{f,g} \Lambda(y)
        &: = 
            \sum\limits_{i=1}^d 
                f^{(i)} (y)
                \tfrac{\partial \Lambda}
                      {\partial x^{(i)}}(y)
            +
            \tfrac{1}{2}
            \sum\limits_{i,l=1}^d 
            \sum\limits_{j=1}^m
                g^{(ij)} (y)
                g^{(lj)} (y)
                \tfrac{\partial^2 \Lambda}
                      {\partial x^{(i)}x^{(l)}}(y),
    \end{aligned}
\end{equation}
where ${\partial \Lambda}/{\partial x^{(i)}}:\mathbb{R}^d
    \rightarrow
    \mathbb{R}^{z_1 \times z_2}$ and 
${\partial^2 \Lambda}/{\partial x^{(i)}x^{(l)}}:\mathbb{R}^d
    \rightarrow
    \mathbb{R}^{z_1 \times z_2}$ represent the first and second order partial derivatives of every element of $\Lambda$. 
For a random variable $\xi: \Omega \rightarrow \mathbb{R}^{z_1 \times z_2}$, we denote $\|\xi\|_{L^r(\Omega;\mathbb{R}^{z_1 \times z_2})}:= (\mathbb{E}[\|\xi\|^r])^{1/r},r>0$.
We said $\Lambda \in \mathcal{C}^1_{\mathcal{P}}(\mathbb{R}^d,\mathbb{R}^{z_1 \times z_2})$ with constants $K_{\Lambda}, c_{\Lambda}$,  if $\Lambda \in  C\big(\mathbb{R}^d,\mathbb{R}^{z_1 \times z_2}\big)$ and 
    \begin{equation}
        \| \Lambda(x)- \Lambda(y))\| \leq K_{\Lambda}(1+|x|+|y|)^{ c_{\Lambda}}|x-y|.  
    \end{equation}
To numerically approximate \eqref{eq:typical_sde}, a uniform mesh with the uniform stepsize $h=\tfrac{T}{N}$, $N \in  \mathbb{Z}^+$, is constructed as follows:
    \begin{equation}\label{eq:uniform_mesh}
        0 = t_{0} < t_{1}< \cdots <t_{N}=T.
    \end{equation} 
For given $N\in \mathbb{Z}^+$, we also define 
    \[
    \lfloor t \rfloor:=\sup_{n=0,...,N}\{t_n:t_n\leq t\},\  t\in [0,T]
    \]
and let a general continuous approximation process defined by
\begin{equation} \label{eq:approxi_sde}
      \left\{ 
        \begin{array}{l}
        \displaystyle{
        Y_t-Y_0=\int_0^t a(s){\rm d}s +\int_0^t b(s){\rm d}W_s,\ t\in [0,T]},\\
        Y_0 = \xi_{Y}.
        \end{array} \right.
\end{equation} 

\section{A new family of stopped increment-tamed methods}
For $t\in[t_k,t_{k+1}],k=0,1,...,N-1$, we propose 
a family of stopped increment-tamed methods starting from $Y_0=\xi_X$ as follows:
\begin{equation}\label{eq:new_stop_tamed_type_methods}
    \begin{aligned}
        &Y_t-Y_{t_k} 
        =\mathbbm{1}_{|Y_{t_k}|\leq \Phi(h)} 
         \Big[
            \Upsilon_h
            \Big(
                \int_{t_k}^t
                F(s,Y_{t_k},W_s)
                {\rm d} s
              +
                \int_{t_k}^t
                G(s,Y_{t_k},W_s)
                {\rm d}W_s
            \Big)
        \Big],
    \end{aligned}
\end{equation}
where
$\Phi:(0,T]\rightarrow (0,\infty)$ 
and  
$\Upsilon_h:\mathbb{R}^d\rightarrow \mathbb{R}^d$  
are two measurable functions given by
$$
    \Phi(x) := \gamma_1 \exp(\gamma_2|\ln(x)|^{\gamma_3}),\gamma_1,\gamma_2>0,0<\gamma_3<1
$$
and
$$
    \Upsilon_h(x) := x\big(1+|x|^{\delta}{h^{-\theta}}\big)^{-1},
    \theta,\delta>0,
$$
respectively. 
Besides, for $t_k \leq s \leq t \leq t_{k+1}$, $F(s,Y_{t_k},W_s)$ and $G(s,Y_{t_k},W_s)$ are 
two adapted processes with respect to the filtration $\{\mathcal{F}_t\}_{t\geq 0}$ and 
rely on $s,Y_{t_k},W_s$. 
Different choices of $F(\cdot,\cdot,\cdot)$ and $G(\cdot,\cdot,\cdot)$ will result in a time-continuous version of different numerical schemes. Typical choices are listed below.
%\newline

%\noindent\textbf{Stopped increment-tamed Euler method:}
By choosing 
$$
    F_1(s,Y_{t_k},W_s)=f(Y_{t_k}),\ G_1(s,Y_{t_k},W_s)=g(Y_{t_k}),
$$ 
one obtains the \textbf{stopped increment-tamed Euler method:}
\begin{equation}\label{eq:stop_tamed_EM_method}
    \begin{aligned}
        &Y_t-Y_{t_k} 
        =\mathbbm{1}_{|Y_{t_k}|\leq \Phi(h)} 
         \Big[
            \Upsilon_h
            \Big(
                \int_{t_k}^t
                f(Y_{t_k})
                {\rm d} s
              +
                \int_{t_k}^t
                g(Y_{t_k})
                {\rm d}W_s
            \Big)
        \Big].
    \end{aligned}
\end{equation}

%\noindent
By setting 
$$
    F_2(s,Y_{t_k},W_s)=f(Y_{t_k}),
    G_2(s,Y_{t_k},W_s)
    = 
        g(Y_{t_k})
        +
        \sum_{j=1}^{m}  
        \mathscr{L}_g^j g(Y_{t_k}) 
        \big(W^{(j)}_s-W^{(j)}_{t_k}\big)
$$ 
one gets the \textbf{stopped increment-tamed Milstein method:}
\begin{equation}\label{eq:stop_tamed_milstein_method}
    \begin{aligned}
        &Y_t-Y_{t_k} 
        =\mathbbm{1}_{|Y_{t_k}|\leq \Phi(h)} 
         \Big[
            \Upsilon_h
            \Big(
                \int_{t_k}^t
                f(Y_{t_k})
                {\rm d} s
              +
                \int_{t_k}^t
                \Big(
                    g(Y_{t_k})
                    +
                    \sum\limits_{j=1}^{m}  
                    \mathscr{L}_g^j g(Y_{t_k}) 
                    \big(W^{(j)}_s-W^{(j)}_{t_k}\big)
                \Big)
                    {\rm d}W_s
            \Big)
        \Big].
    \end{aligned}
\end{equation}

%\noindent
Further, taking
$$
    F_3(s,Y_{t_k},W_s)=f(Y_{t_k})
                        +
                            \sum_{j=1}^{m}
                            \mathscr{L}_g^j f(Y_{t_k})
                            \big(W^{(j)}_s-W^{(j)}_{t_k}\big)
                        +
                        \mathcal{A}_{f,g} f(Y_{t_k}) (s-t_k)
$$
and
\begin{equation*}
    \begin{aligned}
    G_3(s,Y_{t_k},W_s)
    =  
            &g(Y_{t_k})
            +   
            \sum\limits_{j=1}^{m}  
            \mathscr{L}_g^j g(Y_{t_k}) 
            \big(W^{(j)}_s-W^{(j)}_{t_k}\big)
            +
            \mathcal{A}_{f,g} g(Y_{t_k}) (s-t_k)
    \\
            &
            +
            \int_{t_k}^s
            \sum\limits_{j_1,j_2=1}^{m}   
            \mathscr{L}_g^{j_2}
            \mathscr{L}_g^{j_1}
            g(Y_{t_k}) 
            \big(W^{(j_2)}_r-W^{(j_2)}_{t_k}\big)
            {\rm d}W^{(j_1)}_r
    \end{aligned}
\end{equation*}
gives a \textbf{stopped increment-tamed order $1.5$ method:} 
\begin{equation}\label{eq:stop_tamed_1.5_order_method}
    \begin{aligned}
        &Y_t-Y_{t_k} 
        \\
        &=\mathbbm{1}_{|Y_{t_k}|\leq \Phi(h)} 
         \bigg[
            \Upsilon_h
            \bigg(
                \int_{t_k}^t
                \Big(
                    f(Y_{t_k})
                    +
                        \sum_{j=1}^{m}
                            \mathscr{L}_g^j f(Y_{t_k})
                            \big(W^{(j)}_s-W^{(j)}_{t_k}\big)
                    +
                        \mathcal{A}_{f,g} f(Y_{t_k}) (s-t_k)
                \Big)
                {\rm d}s
        \\
        & \quad \quad \quad \quad \quad \quad \quad
            +
            \int_{t_k}^t
            \Big(
                g(Y_{t_k})
                +   
                \sum\limits_{j=1}^{m}  
                \mathscr{L}_g^j g(Y_{t_k}) 
                \big(W^{(j)}_s-W^{(j)}_{t_k}\big)
                +
                \mathcal{A}_{f,g} g(Y_{t_k}) (s-t_k)
        \\
        & \quad \quad \quad \quad
        \quad \quad \quad \quad \quad \quad
                +
                \int_{t_k}^s
                \sum\limits_{j_1,j_2=1}^{m}   
                \mathscr{L}_g^{j_2}
                \mathscr{L}_g^{j_1}
                 g(Y_{t_k}) 
                \big(W^{(j_2)}_r-W^{(j_2)}_{t_k}\big)
                {\rm d}W^{(j_1)}_r
            \Big)
                {\rm d}W_s
            \bigg)
        \bigg].
    \end{aligned}
\end{equation}
For convenience  we introduce the notation of $\tau_N:\Omega \rightarrow \{0,t_1,...,t_N\}$
    \begin{equation}
        \tau_N
        :=
        \inf \big\{ \{T\}
                     \cup 
                     \{     
                         t\in \{ 0,t_1,...,t_N \}:
                         |Y_{t}|> \Phi(h)
                      \}
         \big\}
    \end{equation}
and for $t\in[t_k,t_{k+1}),k=0,1,...,N-1$,
    \begin{equation}
        {Z}_{t} :=
            \int_{t_k}^{t}  F(s,Y_{t_k},W_s) {\rm d}s
            +
            \int_{t_k}^{t} 
                 G(s,Y_{t_k},W_s)
             {\rm d}W_s.
    \end{equation}

\begin{lemma}\label{lem:stop_tamed_method_ito_version}
The equation \eqref{eq:new_stop_tamed_type_methods} can be written as an  It\^{o}'s process  
\begin{equation} \label{eq:stop_tamed_method_ito_version}
      \left\{ 
        \begin{array}{l}
        \displaystyle{
        Y_t-Y_0=\int_0^t a(s){\rm d}s +\int_0^t b(s){\rm d}W_s,\ t\in [0,T]},\\
        Y_0 = \xi_{X},
        \end{array} \right.
    \end{equation}
    where for the Euclidean orthonormal basis $e_1=(1,...,0),...,e_m=(0,...,1)$ and $s \in [t_k,t_{k+1}]$,
\begin{equation}
    \begin{aligned}
        a(s)=\mathbbm{1}_{s< \tau_N}\Big[ \Upsilon_h '(Z_{s})F(s,Y_{t_k},W_s)
        +
        \tfrac{1}{2} \sum\limits_{j=1}^m \Upsilon_h ''(Z_{s})
            \big(
            G(s,Y_{t_k},W_s) e_j,G(s,Y_{t_k},W_s) e_j
            \big)
        \Big]
    \end{aligned}
\end{equation}
and 
\begin{equation}
    \begin{aligned}
        b(s)=\mathbbm{1}_{s< \tau_N} \Upsilon_h '(Z_{s})G(s,Y_{t_k},W_s).
    \end{aligned}
\end{equation}
In addition, for $x,u\in \mathbb{R}^d$,
\begin{equation}\label{eq:Upsilon_one_deri}
        \Upsilon_h '(x)u
        =
        \left\{\begin{array}{ll}u 
            & : x=0 ,
        \\ 
            -\delta h^{-\theta}(1+|x|^{\delta}h^{-\theta})^{-2}|x|^{\delta-2}
            \langle x,u \rangle x
                +
                (1+|x|^{\delta}h^{-\theta})^{-1}u & : x \neq 0,
        \end{array}\right.
\end{equation} 
and
    \begin{equation}\label{eq:Upsilon_two_deri}
       \Upsilon_h ''(x)(u,u)=
       \left\{\begin{array}{ll}0 & : x=0,
       \\ 
            2\delta^2 h^{-2\theta} (1+|x|^{\delta}h^{-\theta})^{-3} |x|^{2\delta-4}|\langle x,u \rangle|^2  x & 
       \\
            \quad 
            +(-\delta) (\delta-2) h^{-\theta} (1+|x|^{\delta}h^{-\theta})^{-2}|x|^{\delta-4}|\langle x,u \rangle|^2  x
       \\ 
            \quad
            +(-\delta)h^{-\theta} (1+|x|^{\delta}h^{-\theta})^{-2} |x|^{\delta-2}
            \big[ 
                u^2x+2\langle x,u \rangle u
            \big]
                 & : x \neq 0.\end{array}\right.
    \end{equation}

\end{lemma}
\textbf{Proof:}
For $t\in[t_k,t_{k+1}], k=0,1,...,N-1$, we denote 
    \begin{equation}
        \tilde{Z}_{k,t}
        :=
            \int_{t_k}^{t} F(s,Y_{t_k},W_s) {\rm d}s
        +
            \int_{t_k}^{t} 
                G(s,Y_{t_k},W_s)
             {\rm d}W_s.
    \end{equation}
The It\^{o} formula then gives
 \begin{equation}
 \begin{aligned}
      \Upsilon_h(\tilde{Z}_{k,t})
      &=\int_{t_k}^{t}
      \Big(
            \Upsilon_h '(\tilde{Z}_{k,s})F(s,Y_{t_k},W_s)
            +
            \tfrac{1}{2}\sum\limits_{j=1}^m \Upsilon_h ''(\tilde{Z}_{k,s})
            \big(G(s,Y_{t_k},W_s) e_j, G(s,Y_{t_k},W_s) e_j\big)
      \Big) 
      {\rm d}s
      \\
            &\quad 
            + \int_{t_k}^{t} 
                \Upsilon_h '(\tilde{Z}_{k,s})
                G(s,Y_{t_k},W_s) 
                {\rm d}W_s
            ,\ t\in[t_k,t_{k+1}].
       \end{aligned}
  \end{equation} 
  Observe that for $x,u,v\in \mathbb{R}^d$ with $x\neq 0$,
  \begin{equation}
    \begin{aligned}
        \Upsilon_h'(x)u
        =
            \sum\limits_{i=1}^d \tfrac{\partial \Upsilon_h}{\partial x^{(i)}}(x) u^{(i)}
        =
            -\delta h^{-\theta}(1+|x|^{\delta}h^{-\theta})^{-2}|x|^{\delta-2}\langle x,u \rangle x+(1+|x|^{\delta}h^{-\theta})^{-1}u,
    \end{aligned}
  \end{equation}
and
 \begin{equation}
    \begin{aligned}
          \Upsilon_h ''(x)(u,v)
          &=
            \sum\limits_{i,l=1}^d 
                \tfrac{\partial^2 \Upsilon_h}{\partial x^{(i)}\partial x^{(l)}}(x) u^{(i)} v^{(l)}
          \\
          &=
            2\delta^2 h^{-2\theta} (1+|x|^{\delta}h^{-\theta})^{-3} |x|^{2\delta-4}\langle x,v \rangle \langle x,u \rangle x
            \\
            &\quad 
                +(-\delta) (\delta-2) h^{-\theta} (1+|x|^{\delta}h^{-\theta})^{-2}|x|^{\delta-4}\langle x,v \rangle \langle x,u \rangle x
            \\
            &\quad
                +(-\delta)h^{-\theta} (1+|x|^{\delta}h^{-\theta})^{-2} |x|^{\delta-2}\big[ 
                    \langle u,v \rangle x+\langle x,u \rangle v+\langle x,v \rangle u
                \big].
    \end{aligned}
  \end{equation}
The proof is thus completed since for any $t\in[t_k,t_{k+1}]$, $\tilde{Z}_{k,t}$ only differs from  $Z_{t}$ at $t=t_{k+1}$. \qed

Next we present some properties of the taming function. Also, some properties of the scheme \eqref{eq:new_stop_tamed_type_methods} are given below, which is useful to prove the exponential integrability properties and strong convergence rate.
\begin{lemma}\label{lem:function_in_schme_property_01}
    For any $\theta, \delta >0 $, it holds
    \begin{equation}
            \sup_{x\in \mathbb{R}} |\Upsilon_h (x)|  \leq  h^{\theta / \delta}.
\end{equation}
Further, for any fixed $v>0$ and $0<h\leq  \exp\big(-\big(\tfrac{\gamma_2}{v}\big)^{1/(1-\gamma_3)}\big)$,
we have
$$
    \gamma_1
    \exp\big({\gamma_2{|\ln (h)|^{\gamma_3}}}\big) \leq \tfrac{\gamma_1}{h^v}.
$$
\end{lemma}
\textbf{Proof:} It is easy to check that 
$$
   |\Upsilon_h (x)|^{\delta}
        \leq
   |x|^{\delta}\big(1+|x|^{\delta}{h^{-\theta}}\big)^{-1} 
        \leq h^\theta.
$$
Moreover, for fixed $v>0$ and $0<\gamma_3<1$, it follows
\begin{equation}
    \begin{aligned}
        \gamma_2 (\ln (\tfrac{1}{h}))^{\gamma_3}
        \leq
         v \ln (\tfrac{1}{h}),
    \end{aligned}
\end{equation}
when $h \leq \exp\big(-\big(\tfrac{\gamma_2}{v}\big)^{1/(1-\gamma_3)}\big)<1$.
This finishes the proof. \qed

\begin{lemma}\label{lem:function_in_schme_property_02}
     For any $x\in \mathbb{R}^d$, it holds that
        \begin{align}
   \|\Upsilon_h'(x)\|_{L(\mathbb{R}^d,\mathbb{R}^d)}
              & \leq 
              \delta |x|^{\delta}{h^{-\theta}}+1;
        \\
              \|\Upsilon_h'(x)-I\|_{L(\mathbb{R}^d;\mathbb{R}^d)} 
              & \leq 
              (\delta+1)|x|^{\delta}{h^{-\theta}};  
        \\
              \sup_{u\in \mathbb{R}^d,|u|\leq 1}|\Upsilon_h''(x)(u,u)| 
              & \leq 
              2\delta^2|x|^{2\delta-1}{h^{-2\theta}}
                +
             (\delta^2+5\delta)|x|^{\delta-1}{h^{-\theta}}.
        \end{align}
\end{lemma}
\textbf{Proof:} The above assertions can be easily validated by \eqref{eq:Upsilon_one_deri} and \eqref{eq:Upsilon_two_deri}. \qed 

\begin{lemma}\label{lem:stop_tamed_method_property}
    Let $Y_t$ be defined by the scheme \eqref{eq:new_stop_tamed_type_methods} for $t\in[t_k,t_{k+1}],k=0,1,...,N-1$. Then
    \begin{enumerate}
        \item for any $t\in[t_k,t_{k+1}],k=0,1,...,N-1$, it holds
        \begin{equation}
        \begin{aligned}
            |Y_t-Y_{t_k}|
            & 
            \leq 
                h^{\theta/ \delta};
        \end{aligned}
    \end{equation}
    \item for any $v_0 \geq 0,\epsilon_0>0$ and $t\in[t_k,t_{k+1}],k=0,1,...,N-1$, there exists some constant $h_0 = \exp\big(-\big(\tfrac{\gamma_2 v_0}{\epsilon_0}\big)^{1/(1-\gamma_3)}\big) > 0$ such that
    \begin{equation}
        \begin{aligned}
            \mathbbm{1}_{|Y_{t_k}|\leq \Phi(h)}|Y_t|^{v_0}
            \leq 
            (2^{v_0-1}+1)(1+\gamma_1^{v_0})h^{-\epsilon_0}, \ h\in (0,h_0].
        \end{aligned}
    \end{equation}
    \end{enumerate}
\end{lemma} 
\textbf{Proof:} 
In view of Lemma \ref{lem:function_in_schme_property_01} one derives that  
\begin{equation}
        \begin{aligned}
%            &
            |Y_t-Y_{t_k}|
%        \\
%           &
           =\mathbbm{1}_{|Y_{t_k}|\leq \Phi(h)}
            \Big[
            \Upsilon_h
            \Big(
                \int_{t_k}^t
                        F(s,Y_{t_k},W_s)
                    {\rm d}s
                +
                    \int_{t_k}^t
                        G(s,Y_{t_k},W_s)
                    {\rm d}W_s
            \Big)
            \Big]
%        \\
%            &
            \leq h^{\theta/ \delta}.
        \end{aligned}
    \end{equation}
    Moreover, for fixed $v_0\geq 0,\epsilon_0>0$, applying Lemma \ref{lem:function_in_schme_property_01} ($v=\epsilon_0/v_0$) again yields 
    \begin{equation}
        \begin{aligned}
            \mathbbm{1}_{|Y_{t_k}|\leq \Phi(h)}|Y_t|^{v_0}
            &\leq 
            \mathbbm{1}_{|Y_{t_k}|\leq \Phi(h)}(|Y_t-Y_{t_k}|+|Y_{t_k}|)^{v_0}
            \\
            & \leq (2^{v_0-1}+1)(h^{\theta v_0/\delta }+\tfrac{\gamma_1^{v_0}}{h^{\epsilon_0}})
            \\
            &\leq (2^{v_0-1}+1)(1+\gamma_1^{v_0})h^{-\epsilon_0}.
        \end{aligned}
    \end{equation}
     The desired assertions are thus validated. \qed
    
\section{Exponential integrability properties}
To show the exponential integrability properties of the proposed scheme, we rely on a key lemma, which can be found in \cite[Corollary 2.3]{hutzenthaler2018exponential}.
\begin{lemma}\label{lem:ex_inte_property_lemma}
    Let a uniform mesh be constructed by \eqref{eq:uniform_mesh}, 
    let $T,\alpha,c\in[0,\infty)$ 
    and let $\psi:\mathbb{R}^d \times [0,T] \times \mathbb{R}^m \rightarrow \mathbb{R}^d$ 
    be a measurable function. 
    Assume that $A\in \mathcal{B}(\mathbb{R}^d)$ and $U_0:\mathbb{R}^d \rightarrow [0,\infty), U_1:\mathbb{R}^d \rightarrow \mathbb{R}$ are  measurable functions. 
    Let $Y:[0,T]\times \Omega \rightarrow \mathbb{R}^d$ be an $\{\mathcal{F}_t\}_{t\in[0,T]}$-adapted stochastic process satisfing that for any $t\in [t_k,t_{k+1}], k=0,1,...,N-1$,
    \begin{equation}
        Y_t=
            \mathbbm{1}_{ \mathbb{R}^d \backslash A}(Y_{t_k})Y_{t_k}
            +
            \mathbbm{1}_{A}(Y_{t_k})\psi(Y_{t_k},t-{t_k},W_t-W_{t_k}),
    \end{equation}
Suppose that for any $x\in A$,
    \begin{equation}
        \int_0^T \mathbbm{1}_{A}(Y_{\lfloor r \rfloor}) |U_1(Y_r)|{\rm d}r
        +
        \int_0^h |U_1(\psi(x,r,W_r))|{\rm d}r 
        <\infty,\ a.s..
    \end{equation}
Besides, assume that for any $(s,x)\in(0,h]\times A$ it holds that
\begin{equation}\label{eq:ex_condi_one_step}
        \mathbb{E}
        \Big[ 
            \exp\Big(
                e^{-\alpha s}U_0(\psi(x,s,W_s))
                +
                    \int_0^s e^{-\alpha r}U_1(\psi(x,r,W_r)){\rm d}r 
                \Big)
        \Big] 
    \leq e^{cs+U_0(x)}.
\end{equation}
Then $\{Y_t\}_{t\in[0,T]}$ admits  the following exponential integrability property
\begin{equation}
    \mathbb{E}
    \Big[ 
    \exp\Big(
        e^{-\alpha t}U_0(Y_t)
        +
        \int_0^t \mathbbm{1}_{A}(Y_{\lfloor r \rfloor}) e^{-\alpha r}U_1(Y_r){\rm d}r 
        \Big) 
    \Big] 
    \leq  e^{ct} 
          \mathbb{E}\big[e^{U_0(Y_0)}\big].
\end{equation}
\end{lemma}

% Now we are well-prepared to show the exponential integrability property of the stopped increment-tamed methods based on Lemma \ref{lem:ex_inte_property_lemma}.

In the following, based on Lemma \ref{lem:ex_inte_property_lemma}, we construct a unified framework to validate the exponential integrability properties of a family of stopped increment-tamed methods \eqref{eq:new_stop_tamed_type_methods}.

\begin{theorem}\label{thm:ex_inte_property_of_stop_tamed_method} 
Let 
    $U_0\in {C}^2(\mathbb{R}^d,[0,\infty)) \cap  \mathcal{C}^1_{\mathcal{P}}(\mathbb{R}^d,[0,\infty))  $ with constants $K_{0},c_{0}\geq 0$ 
and let 
$U_1 \in C(\mathbb{R}^d,\mathbb{R})$ satisfy 
$|U_1(x)|\leq K_{1}(1+|x|)^{c_{1}}, K_{1},c_{1}\geq 0$.
Let 
    $f\in \mathcal{C}^1_{\mathcal{P}}(\mathbb{R}^d,\mathbb{R}^d)$
and 
    $g \in \mathcal{C}^1_{\mathcal{P}}(\mathbb{R}^d,\mathbb{R}^{d\times m})$
with constants $K_f,c_f$ and $K_g,c_g$, respectively.
Let a family of stopped increment-tamed methods be defined by \eqref{eq:new_stop_tamed_type_methods} fulfilling 
$$
   \mathbb{E}\left[e^{U_0(\xi_X )}\right]<\infty.
$$
Assume that there exists constant $\alpha \in [0,\infty)$ such that for any $y\in \mathbb{R}^d$,
\begin{equation}\label{eq:ex_inte_property_condition}
    \begin{aligned}
            \mathcal{A}_{f,g} U_0(y)
                +\tfrac{1}{2}|g(y)^{*}(\nabla U_0(y))|^2
                +U_1(y) 
            \leq 
                \alpha U_0(y).
    \end{aligned}
\end{equation}
Further, suppose that 
for any $s\in (0,h], |x|\leq \gamma_1 \exp(\gamma_2|\ln h|^{\gamma_3})$ and $Y_0=x$, there exist some constant $C_a, C_{b}, \varepsilon_a>0, \tfrac{2\theta}{\delta} > \varepsilon_{b}>0$ such that
\begin{equation}\label{eq:ex_inte_property_condi_on_a}
    \begin{aligned}
            \big\|
            |f(x)-a(s)|
            \big \|_{L^4(\Omega;\mathbb{R})} \leq C_a h^{\varepsilon_a}
    \end{aligned}
\end{equation}
and
\begin{equation}\label{eq:ex_inte_property_condi_on_b}
    \begin{aligned}
        \big\|
            \|b(s)\|
        \big \|_{L^8(\Omega;\mathbb{R})} \leq C_b(1+ h^{-\frac{\varepsilon_b}{2}}),
        \big\|
            \|g(x)-b(s)\|
        \big \|_{L^8(\Omega;\mathbb{R})} \leq C_b h^{\varepsilon_b},  
    \end{aligned}
\end{equation}
when $h$ is sufficiently small.
Then it holds that
\begin{equation}\label{eq:ex_inter_property_01}
    \begin{aligned}
    &\sup_{N\in \mathbb{Z}^+}
    \sup_{t\in[0,T]} 
    \mathbb{E}
    \Big[ 
        \exp\Big(
            e^{-\alpha t}U_0(Y_t)+\int_0^{t\wedge \tau_N} e^{-\alpha r}U_1(Y_r){\rm d}r 
            \Big) 
    \Big] 
    \leq  e^{CT} \mathbb{E}\big[
                    e^{U_0(\xi_X)}
                    \big]
    \end{aligned}
\end{equation}
for some constant $C>0$ and 
\begin{equation}\label{eq:ex_inter_property_02}
    \begin{aligned}
    &\limsup_{ N \rightarrow +\infty}
    \sup_{t\in[0,T]}
    \mathbb{E}
    \Big[ 
        \exp\Big(
            e^{-\alpha t}U_0(Y_t)+\int_0^{t\wedge \tau_N} e^{-\alpha r}U_1(Y_r){\rm d}r 
            \Big) 
    \Big]
    \leq  \mathbb{E}\big[
                    e^{U_0(\xi_X)}
                    \big].
    \end{aligned}
\end{equation}

\end{theorem}
\textbf{Proof:}
Firstly, for $t\in (0,h]$ and 
    $|x|\leq \gamma_1 \exp(\gamma_2|\ln h|^{\gamma_3})$, we denote
    \begin{equation}
        \begin{aligned}
            Y^x_t:=
            \psi(x,t,W_t) 
            =x+ 
                \Upsilon_h 
                    \Big( \int_0^t  
                            F(s,x,W_s) {\rm d}s
                            + \int_{0}^t 
                                 G(s,x,W_s)
                              {\rm d}W_s
                    \Big).
        \end{aligned}
    \end{equation}
    To confirm condition \eqref{eq:ex_condi_one_step},
    by It\^{o}'s formula and equation \eqref{eq:ex_inte_property_condition} one obtains that
\begin{align*}
            &\exp\Big(
                    e^{-\alpha t}U_0(Y^x_t)+\int_0^{t} e^{-\alpha r} U_1(Y^x_r) {\rm d} r
                 \Big)
             -e^{U_0(x)}
        \\
            &= \int_{0}^{t}
                 \exp\Big(
                        e^{-\alpha s}U_0(Y^x_s)+\int_0^{s}   e^{-\alpha r} U_1(Y^x_r) {\rm d} r
                    \Big) 
                    e^{-\alpha s} 
                    U_0'(Y^x_s)b(s) 
                {\rm d}W_s
        \\
            &\quad 
            + \int_{0}^{t }  
            \exp\Big(
                    e^{-\alpha s}U_0(Y^x_s)+\int_0^{s}   e^{-\alpha r} U_1(Y^x_r) {\rm d} r
                \Big) 
            e^{-\alpha s} 
            \cdot\Big[
                    U_0'(Y^x_s)a(s)
        \\
            &\quad \quad 
            +\tfrac{1}{2}
                \operatorname{tr}\big(b(s)^*\operatorname{Hess}_x(U_0(Y^x_s))b(s)\big)
            +\tfrac{1}{2e^{\alpha s}}|b^*(s)\nabla U_0(Y^x_s)|^2 
            +U_1(Y^x_s)
            -\alpha U_0(Y^x_s)
            \Big] 
            {\rm d}s
        \\
            &= \int_{0}^{t}
                 \exp\Big(
                        e^{-\alpha s}U_0(Y^x_s)+\int_0^{s}   e^{-\alpha r} U_1(Y^x_r) {\rm d} r
                    \Big) 
                    e^{-\alpha s} 
                    U_0'(Y^x_s)b(s) 
                {\rm d}W_s
        \\
            &\quad 
            + \int_{0}^{t }  
                \exp\Big(
                    e^{-\alpha s}U_0(Y^x_s)+\int_0^{s} e^{-\alpha r} U_1(Y^x_r) {\rm d} r
                \Big) 
                e^{-\alpha s}
                \cdot\Big[U_0'(Y^x_s)f(Y^x_s) 
        \\
            &\quad \quad 
            +\tfrac{1}{2}\operatorname{tr}\big(g(Y^x_s)^*\operatorname{Hess     }_x(U_0(Y^x_s))g(Y^x_s)\big)
            +\tfrac{1}{2e^{\alpha s}}|g^*(Y^x_s)\nabla U_0(Y^x_s)|^2 
            + U_1(Y^x_s)
            -\alpha U_0(Y^x_s)
            \Big] 
            {\rm d}s 
        \\
            &\quad 
            + \int_{0}^{t } \exp\Big(
                e^{-\alpha s}U_0(Y^x_s)+\int_0^{s} e^{-\alpha r} U_1(Y^x_r) {\rm d}r
                \Big) 
                e^{-\alpha s} \cdot\Big[
                    U_0'(Y^x_s)
                    \big(a(s)-f(Y^x_s)\big)
        \\
            &\quad \quad 
            +\tfrac{1}{2}\big[
            \operatorname{tr}\big(b(s)^*\operatorname{Hess}_x(U_0(Y^x_s))b(s)\big)
                -\operatorname{tr}\big(g(Y^x_s)^*\operatorname{Hess}_x(U_0(Y^x_s))g(Y^x_s)
            \big)\big]
        \\
            &\quad \quad 
            +\tfrac{1}{2e^{\alpha s}}\big(
            |b^*(s)\nabla U_0(Y^x_s)|^2-|g^*(Y^x_s)\nabla U_0(Y^x_s)|^2
            \big)
            \Big] {\rm d}s 
        \\
            &\leq \int_{0}^{t}
                 \exp\Big(
                        e^{-\alpha s}U_0(Y^x_s)+\int_0^{s}   e^{-\alpha r} U_1(Y^x_r) {\rm d} r
                    \Big) 
                    e^{-\alpha s} 
                    U_0'(Y^x_s)b(s) 
                {\rm d}W_s
        \\
            &\quad 
            + \int_{0}^{t } \exp\Big(
                e^{-\alpha s}U_0(Y^x_s)+\int_0^{s} e^{-\alpha r} U_1(Y^x_r) {\rm d}r
                \Big) 
                e^{-\alpha s} \cdot\Big[
                    U_0'(Y^x_s)
                    \big(a(s)-f(Y^x_s)\big)
        \\
            &\quad \quad 
            +\tfrac{1}{2}\big[
            \operatorname{tr}\big(b(s)^*\operatorname{Hess}_x(U_0(Y^x_s))b(s)\big)
                -\operatorname{tr}\big(g(Y^x_s)^*\operatorname{Hess}_x(U_0(Y^x_s))g(Y^x_s)
            \big)\big]
        \\
            &\quad \quad  
            +\tfrac{1}{2e^{\alpha s}}\big(
            |b^*(s)\nabla U_0(Y^x_s)|^2-|g^*(Y^x_s)\nabla U_0(Y^x_s)|^2
            \big)
            \Big] {\rm d}s.
            \stepcounter{equation}
            \tag{\theequation}
      \end{align*} 
Using the localization technique of stopping time and 
taking expectation of both side yields  
    \begin{align*}
            &\mathbb E \Big[ \exp\Big(e^{-\alpha t}U_0(Y^x_t)+\int_0^{t} e^{-\alpha r} U_1(Y^x_r) {\rm d} r\Big )\Big]
         - e^{U_0(x)}
        \\
            &\leq \mathbb{E} \bigg[
            \int_{0}^{t } \exp\Big(
                e^{-\alpha s}U_0(Y^x_s)+\int_0^{s} e^{-\alpha r} U_1(Y^x_r) {\rm d}r
                \Big) 
                e^{-\alpha s} \cdot\Big|
                    U_0'(Y^x_s)
                    \big(a(s)-f(Y^x_s)\big)
        \\
            &\quad \quad \quad 
            +\tfrac{1}{2}\big[
            \operatorname{tr}\big(b(s)^*\operatorname{Hess}_x(U_0(Y^x_s))b(s)\big)
                -\operatorname{tr}\big(g(Y^x_s)^*\operatorname{Hess}_x(U_0(Y^x_s))g(Y^x_s)
            \big)\big]
        \\
            &\quad \quad \quad 
            +\tfrac{1}{2e^{\alpha s}}\big(
            |b^*(s)\nabla U_0(Y^x_s)|^2-|g^*(Y^x_s)\nabla U_0(Y^x_s)|^2
            \big)
            \Big| {\rm d}s 
            \bigg]
        \\
            &\leq  e^{U_0(x)} 
            \bigg[
            \int_{0}^{t}  
            \bigg(
            \underbrace{
                \Big\| 
                    \exp\Big(
                    e^{-\alpha s}U_0(Y^x_s)
                    -{U_0(x)}
                    +\int_{0}^{s} e^{-\alpha r} U_1(Y^x_r) {\rm d} r
                    \Big)
                \Big\|_{L^2(\Omega;\mathbb{R})} 
            }_{=:J_0}
        \\
            &\quad \quad \quad \quad 
            \cdot \Big[
            \underbrace{
                \big\|
                    U_0'(Y^x_s)(a(s)-f(Y^x_s))
                \big\|_{L^2(\Omega;\mathbb{R})} 
            }_{=:J_1}
            +
            \underbrace{
            \tfrac{1}{2}
            \big\| 
                |b^*(s)\nabla U_0(Y^x_s)|^2
                -|g^*(Y^x_s)\nabla U_0(Y^x_s)|^2
            \big\|_{L^2(\Omega;\mathbb{R})} 
            }_{=:J_2}
        \\
            &\quad \quad \quad \quad \quad +
            \underbrace{
            \tfrac{1}{2}
            \big\|
                \operatorname{tr}\big(b(s)^*\operatorname{Hess}_x(U_0(Y^x_s))b(s)\big)
                -
                \operatorname{tr}\big(g(Y^x_s)^*\operatorname{Hess}_x(U_0(Y^x_s))g(Y^x_s)
            \big\|_{L^2(\Omega;\mathbb{R})} 
            }_{=:J_3} 
            \Big] 
            \bigg) 
            {\rm d}s
             \bigg]. 
    \stepcounter{equation}\tag{\theequation}
        \label{eq:estimate_stopped_tamed_method_J0_to_J3}
    \end{align*}
One needs to estimate $J_0,J_1,J_2,J_3$ term by term. 
For $J_0$, by the conditions on $U_0$ and $U_1$, $e^{-x}-1 \leq x \ (x\geq 0)$, Lemma \ref{lem:function_in_schme_property_01} and Lemma \ref{lem:stop_tamed_method_property}, there exists some constant $h_0>0$ such that for any $0< h = \tfrac{T}{N} \leq h_0$
\begin{equation} \label{eq:process_estimate_stop_tamed_method_J0_01}
    \begin{aligned}
    J_0& = 
        \Big\| 
        \exp\Big(
        e^{-\alpha s}(U_0(Y^x_s)-U_0(x))
        +
        (e^{-\alpha s}-1)U_0(x)
        +
        \int_{0}^{s} e^{-\alpha r} U_1(Y^x_r) {\rm d} r
        \Big)
        \Big\|_{L^2(\Omega;\mathbb{R})} 
    \\
        &\leq 
        \Big\| 
        \exp\Big(
        K_0(1+|Y^x_s|+|x|)^{c_{0}}h^{\frac{\theta}{\delta}}
        +
        {\alpha h}\big(K_0(1+|x|)^{c_{0}+1} 
                    +|U_0(0)|\big)
    \\
        &\quad\quad \quad\quad  + 
        \int_{0}^{s} 
            K_1(1+|Y^x_r|)^{c_{1}}
             {\rm d} r
        \Big) 
        \Big\|_{L^2(\Omega;\mathbb{R})} 
    \\
        &\leq 
        \exp\Big(
        K_0 (1+3^{c_0-1})
        \big(
                1+(2^{c_0-1}+1)(1+\gamma_1^{c_0})h^{-\frac{\theta}{2\delta}} +
                \gamma_1 h^{-\frac{\theta}{2\delta}}
        \big) h^{\frac{\theta}{\delta}}
    \\
        &\quad\quad \quad\quad 
        +
        {\alpha h}
        \big(2^{c_0}K_0
        (
            1+
            \gamma_1 h^{-\frac{1}{2}}
        )+|U_0(0)|\big)
    \\
        &\quad\quad \quad\quad 
        + 
        h K_1 
        (2^{c_1-1}+1)
            \big(
                    1
                +
                    (2^{c_1-1}+1)
                    (1+\gamma_1^{c_1})
                    h^{-\tfrac{1}{2}} 
            \big)
        \Big). 
    \end{aligned}
\end{equation}
Meanwhile, it is clear that there exists a constant $C_0>0$ fulfilling
\begin{equation}\label{eq:process_estimate_stop_tamed_method_J0_big_h}
    \max_{h_0 < h \leq T} J_0 \leq C_0.
\end{equation} 
One then concludes that
\begin{equation}\label{eq:result_estimate_stop_tamed_method_J0}
    J_0 \leq C_0(h),
\end{equation} 
where $C_0(h)>0$ is bounded with respect to $h$ and satisfies 
$$
    \lim_{h\rightarrow 0^+} C_0(h)=1.
$$
For $J_1$, noticing $U_0 \in \mathcal{C}^1_{\mathcal{P}}(\mathbb{R}^d,[0,\infty)),f\in \mathcal{C}^1_{\mathcal{P}}(\mathbb{R}^d,\mathbb{R}^d)$ and using Lemma \ref{lem:stop_tamed_method_property}, the H\"{o}lder inequality and \eqref{eq:ex_inte_property_condi_on_a} show
\begin{equation}\label{eq:process_estimate_stop_tamed_method_J1_01}
    \begin{aligned}
        &
            \|U_0'(Y^x_s)\|_{L(\mathbb{R}^d,\mathbb{R}^d)}
            |f(Y^x_s)-f(x)| 
        \\
            &\leq 
                K_0(1+2|Y^x_s|)^{c_{0}}K_f(1+|Y^x_s|+|x|)^{c_f}
                h^{\frac{\theta}{\delta}}
    \end{aligned}
\end{equation}
and
\begin{equation}\label{eq:process_estimate_stop_tamed_method_J1_02}
    \begin{aligned}
            & 
            \big\| 
            \|U_0'(Y^x_s)\|_{L(\mathbb{R}^d,\mathbb{R}^d)}|f(x)-a(s)|
            \big\|_{L^2(\Omega;\mathbb{R})} 
        \\
            &\leq 
              C_a 
              h^{\varepsilon_a} K_0
              \big\| 
                    (1+2|Y^x_s|)^{c_{0}}
                 \big\|_{L^4(\Omega;\mathbb{R})}, 
    \end{aligned}
\end{equation}
when $h$ is sufficient small.
By repeating a similar argument in the last inequality of
\eqref{eq:process_estimate_stop_tamed_method_J0_01} 
and
\eqref{eq:process_estimate_stop_tamed_method_J0_big_h} one acquires
\begin{equation}\label{eq:result_estimate_stop_tamed_method_J1}
    \begin{aligned}
        &J_1 \leq C_1(h). 
    \end{aligned}
\end{equation}
Here $C_1(h)>0$ is bounded with respect to $h$ and satisfies 
$$
    \lim_{h\rightarrow 0^+} C_1(h)=0.
$$
With regard to $J_2$, we use $U_0 \in \mathcal{C}^1_{\mathcal{P}}(\mathbb{R}^d,[0,\infty))$, the H\"{o}lder inequality, Lemma \ref{lem:stop_tamed_method_property} and \eqref{eq:ex_inte_property_condi_on_b} to derive
\begin{equation}
    \begin{aligned}
        J_2 
        &=\tfrac{1}{2}
        \big\| 
                \langle 
                    b(s),\nabla U_0(Y^x_s)U'(Y^x_s)b(s) 
                \rangle_{HS}
              -
                \langle 
                    g(Y^x_s),\nabla U_0(Y^x_s)U'(Y^x_s)g(Y^x_s) 
                \rangle_{HS} 
        \big\|_{L^2(\Omega;\mathbb{R})} 
        \\
        &= 
        \tfrac{1}{2}
        \big\| 
                \langle 
                    b(s)+g(Y^x_s),
                    \nabla U_0(Y^x_s)U_0'(Y^x_s)(b(s)-g(Y^x_s)) 
                \rangle_{HS} 
        \big\|_{L^2(\Omega;\mathbb{R})} 
        \\
        &\leq 
        \tfrac{1}{2}
        \big\| 
                \|
                    \nabla U_0(Y^x_s)U_0'(Y^x_s)
                \|_{L(\mathbb{R}^d,\mathbb{R}^d)} 
        \big\|_{L^4(\Omega;\mathbb{R})} 
        \big\| 
                \| b(s)+g(Y^x_s)\|
        \big\|_{L^8(\Omega;\mathbb{R})} 
        \big\| 
                \| b(s)-g(Y^x_s)\|
        \big\|_{L^8(\Omega;\mathbb{R})} 
        \\
        &\leq 
        C_{u_0}
        \Big(
            1+
            \big\| |Y^x_s|^{2c_{0}} \big\|_{L^4(\Omega;\mathbb{R})}
        \Big)
        \cdot
        \Big[ 
        C_b(1+h^{-\frac{\varepsilon_b}{2}})
            +
        \big\| 
            K_g(1+|Y^x_s|)^{c_{g}+1}+\|g(0)\|
        \big\|_{L^8(\Omega;\mathbb{R})}
        \Big]
        \\
        &
        \quad 
        \cdot 
        \Big[
         \big\| 
                \| b(s)-g(x)\|
        \big\|_{L^8(\Omega;\mathbb{R})}  
        +
        K_g 
        \big\| 
             (1+|x|+|Y^x_s|)^{c_g}
        \big\|_{L^8(\Omega;\mathbb{R})}
        h^{\frac{\theta}{\delta}}
        \Big], 
    \end{aligned}
\end{equation}
where $C_{u_0}$ is some constant depending on $U_0$.
Using \eqref{eq:ex_inte_property_condi_on_b} and the same arguments as used in the last inequality of
\eqref{eq:process_estimate_stop_tamed_method_J0_01}
shows 
\begin{equation}\label{eq:result_estimate_stop_tamed_method_J2}
    \begin{aligned}
        &J_2 \leq C_2(h),
    \end{aligned}
\end{equation}
where $C_2(h)>0$ is bounded with respect to $h$ and satisfies 
$$
    \lim_{h\rightarrow 0^+} C_2(h)=0.
$$
The estimate of $J_3$ is analogous to that of $J_2$. 
Observe that  
\begin{equation}
    \begin{aligned}
        J_3 
        &=
        \tfrac{1}{2}
            \big\|
                \langle 
                b(s),\operatorname{Hess}_x(U_0(Y^x_s))b(s)
                \rangle_{HS}
            -
              \langle 
                g(Y^x_s),\operatorname{Hess}_x(U_0(Y^x_s))g(Y^x_s)
              \rangle_{HS}
            \big\|_{L^2(\Omega;\mathbb{R})} 
        \\
        &\leq 
        \tfrac{1}{2}
        \big\| 
            \| 
                b(s)+g(Y^x_s)
            \|
            \|
                \operatorname{Hess}_x(U_0(Y^x_s))
            \|_{L(\mathbb{R}^d \times \mathbb{R}^d,\mathbb{R}^d)} 
            \|
                b(s)-g(Y^x_s)
            \|
        \big\|_{L^2(\Omega;\mathbb{R})}.
    \end{aligned}
\end{equation}
Hence there exists some bounded function $C_3(h)>0$ such that 
$\lim_{h\rightarrow 0^+} C_3(h)=0$ and
\begin{equation}\label{eq:result_estimate_stop_tamed_method_J3}
    \begin{aligned} 
        &J_3 \leq C_3(h).
    \end{aligned}
\end{equation}
Therefore, combining 
\eqref{eq:estimate_stopped_tamed_method_J0_to_J3}, 
\eqref{eq:result_estimate_stop_tamed_method_J0},
\eqref{eq:result_estimate_stop_tamed_method_J1},
\eqref{eq:result_estimate_stop_tamed_method_J2} 
with
\eqref{eq:result_estimate_stop_tamed_method_J3}
gives
\begin{equation}
    \begin{aligned}
            &\mathbb E 
            \Big[ 
                \exp\Big(
                    e^{-\alpha t}U_0(Y^x_t)
                        +\int_0^{t} 
                            e^{-\alpha r} U_1(Y^x_r)
                         {\rm d} r
                \Big)
            \Big]
        \\
            &\leq e^{U_0(x)}
            \Big(
            1+C_0(h)
                \big(
                    C_1(h)
                    +
                    C_2(h)
                    +
                    C_3(h)
                \big)
                t
            \Big)
        \\
            &\leq 
            e^{U_0(x)+C_0(h)
                (
                    C_1(h)
                    +
                    C_2(h)
                    +
                    C_3(h)
                )t},
    \end{aligned}
\end{equation}
which confirms equation \eqref{eq:ex_condi_one_step} and 
finally completes the proof  due to  Lemma \ref{lem:ex_inte_property_lemma}. \qed

\begin{corollary}\label{cor:numerical-exponential-moment}
Let 
    $U_0\in {C}^2(\mathbb{R}^d,[0,\infty)) \cap \mathcal{C}^1_{\mathcal{P}}(\mathbb{R}^d,[0,\infty)) $ with constants $K_{0},c_{0}\geq 0$ 
and let 
$U_1 \in C(\mathbb{R}^d,\mathbb{R})$ satisfy 
$|U_1(x)|\leq K_{1}(1+|x|)^{c_{1}}, K_{1},c_{1}\geq 0$.
Let 
    $f\in \mathcal{C}^1_{\mathcal{P}}(\mathbb{R}^d,\mathbb{R}^d)$
and 
    $g \in \mathcal{C}^1_{\mathcal{P}}(\mathbb{R}^d,\mathbb{R}^{d\times m})$
with constants $K_f,c_f$ and $K_g,c_g$, respectively.
Let a family of stopped increment-tamed methods be defined by \eqref{eq:new_stop_tamed_type_methods} fulfilling 
$$
   \mathbb{E}\left[e^{U_0(\xi_X )}\right]<\infty
$$
and $\delta-2\theta >1$.
Assume that there exists constant $\alpha \in [0,\infty)$ such that for any $y\in \mathbb{R}^d$,
\begin{equation}
    \begin{aligned}
            \mathcal{A}_{f,g} U_0(y)
                +\tfrac{1}{2}|g(y)^{*}(\nabla U_0(y))|^2
                +U_1(y) 
            \leq 
                \alpha U_0(y).
    \end{aligned}
\end{equation}
Then the stopped increment-tamed Euler method \eqref{eq:stop_tamed_EM_method},
stopped increment-tamed Milstein method \eqref{eq:stop_tamed_milstein_method}
and
stopped increment-tamed order $1.5$ method \eqref{eq:stop_tamed_1.5_order_method}
admit the exponential integrability properties 
\eqref{eq:ex_inter_property_01} and 
\eqref{eq:ex_inter_property_02}.
\end{corollary}
\textbf{Proof:}
According to Theorem \ref{thm:ex_inte_property_of_stop_tamed_method},  one only needs to validate \eqref{eq:ex_inte_property_condi_on_a} and \eqref{eq:ex_inte_property_condi_on_b}. To do this, let $s \in (0,h]$ and 
$|x|\leq \gamma_1 \exp(\gamma_2|\ln h|^{\gamma_3})$.
For the stopped increment-tamed Euler method \eqref{eq:stop_tamed_EM_method},
the Lemma \ref{lem:stop_tamed_method_ito_version} and  Lemma \ref{lem:function_in_schme_property_02} show
\begin{equation}
    \begin{aligned}
                \big\||f(x)-a(s)|\big\|_{L^4(\Omega;\mathbb{R})}
            &=\big\|
                |
                (\Upsilon_h'(Z_s)-I)f(x)
               + 
               \tfrac{1}{2}\sum\limits_{j=1}^m \Upsilon_h ''(Z_{s})
                \big(
                    g(x)e_j,g(x)e_j
                \big)
                |
             \big\|_{L^4(\Omega;\mathbb{R})}
             \\
             &\leq 
                (\delta +1)
                h^{-\theta} 
                    |f(x)|
                \big\|
                    |Z_s|^{\delta}
                \big\|_{L^4(\Omega;\mathbb{R})}
                +
                \tfrac{m}{2} \cdot
            \\
            & \quad \ 
              \big(
                  2\delta^2 {h^{-2\theta}}
                        \big\|
                            |Z_s|^{2\delta-1}
                        \big\|_{L^4(\Omega;\mathbb{R})}
                  +
                  (\delta^2+5\delta){h^{-\theta}}
                  \big\|
                  |Z_s|^{\delta-1}
                  \big\|_{L^4(\Omega;\mathbb{R})}
              \big)
              \|g(x)\|^2,
    \end{aligned}
\end{equation}
where $Z_s=f(x)s+g(x)W_s$.
For any $\lambda_1 >0$, using the Burkholder-Davis-Gundy inequality gives
\begin{equation}
    \begin{aligned}
        \mathbb{E} \big[ |{Z}_s|^{\lambda_1} \big]
        &
        \leq 
        (2^{{\lambda_1}-1}+1)\Big(
            |f(x)|^{\lambda_1} h^{\lambda_1}
        +
            \|g(x)\|^{\lambda_1}
            \mathbb{E} \big[|W_s|^{\lambda_1}\big]
        \Big)
    \\
        &
        \leq 
        C_{\lambda_1}h^{\frac{\lambda_1}{2}} 
        \Big(
            |f(x)|^{\lambda_1} h^{\frac{\lambda_1}{2}} 
            +
            |g(x)|^{\lambda_1}
            \Big).
    \end{aligned}
\end{equation}
Since $\delta - 2\theta >1$, we then  conclude that 
there exist some constant $C_a, \varepsilon_a>0$ such that
\begin{equation}
    \begin{aligned}
            \big\|
            |f(x)-a(s)|
            \big \|_{L^4(\Omega;\mathbb{R})} \leq C_a h^{\varepsilon_a}.
    \end{aligned}
\end{equation}
Meanwhile, recall that $b(s)=\Upsilon_h'(Z_s)g(x)$. Then
\begin{equation}
    \begin{aligned}
        \big \|
            \|g(x)-b(s)\|
        \big \|_{L^8(\Omega;\mathbb{R})}
        =
        \big \|
            \|
                \Upsilon_h'(Z_s)-I)g(x)
            \|
        \big \|_{L^8(\Omega;\mathbb{R})}.
    \end{aligned}
\end{equation}
By Lemma \ref{lem:function_in_schme_property_02}, it is obvious that 
there exist some constant $C_b, \varepsilon_b>0$ such that for small $h$,
\begin{equation}
    \begin{aligned}
        \big\|
            \|b(s)\|
        \big \|_{L^8(\Omega;\mathbb{R})} \leq C_b(1+ h^{-\frac{\varepsilon_b}{2}}),\
        \big\|
            \|g(x)-b(s)\|
        \big \|_{L^8(\Omega;\mathbb{R})} \leq C_b h^{\varepsilon_b}.  
    \end{aligned}
\end{equation}
For the stopped increment-tamed Milstein method \eqref{eq:stop_tamed_milstein_method},
we rely on Lemma \ref{lem:stop_tamed_method_ito_version} and  Lemma \ref{lem:function_in_schme_property_02} to derive
\begin{equation}
    \begin{aligned}
            &\big\||f(x)-a(s)|\big\|_{L^4(\Omega;\mathbb{R})}
            \\
            &=\big\|
                |
                (\Upsilon_h'({Z}_s)-I)f(x)
               + 
               \tfrac{1}{2}\sum\limits_{j=1}^m \Upsilon_h ''({Z}_{s})
                \big(
                    G_2(x,W_s)e_j,G_2(x,W_s)e_j
                \big)
                |
             \big\|_{L^4(\Omega;\mathbb{R})}
             \\
             &\leq 
                C_{\delta}
                h^{-\theta} 
                    |f(x)|
                \big\|
                    |{Z}_s|^{\delta}
                \big\|_{L^8(\Omega;\mathbb{R})}
                +
                \tfrac{m}{2} \cdot
            \\
            & \quad \quad 
              C_{\delta}\big(
                   {h^{-2\theta}}
                        \big\|
                            |{Z}_s|^{2\delta-1}
                        \big\|_{L^8(\Omega;\mathbb{R})}
                  +
                  {h^{-\theta}}
                  \big\|
                  |{Z}_s|^{\delta-1}
                  \big\|_{L^8(\Omega;\mathbb{R})}
              \big)
              \big\|
                    \| G_2(x,W_s) \|^2
              \big\|_{L^8(\Omega;\mathbb{R})},
    \end{aligned}
\end{equation}
where 
$$
    G_2(x,W_s) = 
            g(x)+\sum_{j=1}^{m}  
            \mathscr{L}_g^j g(x) 
            \big(W^{(j)}_s-W^{(j)}_{t_k}\big),
            \
            {Z}_s= 
            f(x)s+\int_0^s G_2(x,W_r){\rm {d}}W_r.
$$
For any $\lambda_0 \geq 1,\lambda_1 > 0,s\in(0,h]$, one can use the Burkholder-Davis-Gundy inequality to obtain
\begin{equation}\label{eq:estimate_milstein_b_s}
    \begin{aligned}
        &\mathbb{E}  
            \big[
                    \|G_2(x,W_s)\|^{\lambda_0}
            \big]
        \\
        &\leq 
            2^{{\lambda_0}-1} 
            \Big(
                \big(K_g(1+|x|)^{c_g+1}+|g(0)|\big)^{\lambda_0}
                +
                m^{{\lambda_0}-1}
                \sum_{j=1}^m 
                \|\mathscr{L}^j g(x)\|^{\lambda_0}
                \mathbb{E} 
                \big[
                    |W_s^{(j)}|^{{\lambda_0}}
                \big]
            \Big)
        \\
        &\leq 
            2^{{\lambda_0}-1} 
            \Big(
                \big(K_g(1+|x|)^{c_g+1}+|g(0)|\big)^{\lambda_0}
                +
                C_g m^{{\lambda_0}}
                (1+|x|)^{(2c_g+1){\lambda_0}}
                C_{\lambda_0}h^{{\lambda_0}/2}
            \Big)
    \end{aligned}
\end{equation}
and 
\begin{equation}
    \begin{aligned}
        \mathbb{E} \big[ |{Z}_s|^{\lambda_1} \big]
        &
        \leq 
        (2^{{\lambda_1}-1}+1)\Big(
            |f(x)|^{\lambda_1} h^{\lambda_1}
        +
            C_{\lambda_1}\mathbb{E} 
            \Big[
                    \Big|
                        \int_{0}^{s} 
                            \|
                                G_2(x,W_r) 
                            \|^2
                            {\rm d}r
                    \Big|^{\frac{\lambda_1}{2}}
            \Big]
        \Big)
    \\
        &
        \leq 
        C_{\lambda_1}h^{\frac{\lambda_1}{2}} 
        \Big(
            |f(x)|^{\lambda_1} h^{\frac{\lambda_1}{2}} 
            +
                \mathbbm{1}_{\lambda_1 <2}
                \big(\mathbb{E} [
                            \|
                                G_2(x,W_r) 
                            \|^2
                           ]
                \big)^{\frac{\lambda_1}{2}} 
                +
                \mathbbm{1}_{2 \leq \lambda_1}
                \mathbb{E} [
                            \|
                                G_2(x,W_r) 
                            \|^{\lambda_1}
                           ]
            \Big),
    \end{aligned}
\end{equation}
which validates \eqref{eq:ex_inte_property_condi_on_a}. 
The first assertion of \eqref{eq:ex_inte_property_condi_on_b} can be easily deduced by \eqref{eq:estimate_milstein_b_s} and Lemma \ref{lem:function_in_schme_property_02}. 
Moreover, one can infer
\begin{equation}
    \begin{aligned}
        &\big\|
            \|g(x)-b(s)\|
        \big \|_{L^8(\Omega;\mathbb{R})}
    \\
        &=  
        \big\|
            \big(
                \Upsilon_h' ({Z}_s)-I
            \big)
            g(x)
            +
            \Upsilon_h' ({Z}_s)
            \sum_{j=1}^{m}  
            \mathscr{L}_g^j g(x) 
            \big(W^{(j)}_s-W^{(j)}_{t_k}\big)
            \|
        \big \|_{L^8(\Omega;\mathbb{R})}
    \\
        &\leq 
        C_{\delta}
            h^{-\theta}  
            \big\| 
                 |{Z}_s|^{\delta} 
            \big\|_{L^8(\Omega;\mathbb{R})}
             \| g(x)\|
            +   
            C_g
             C_{\lambda_0}
            m
            \big(
            1+ 
            \delta h^{-\theta}
            \big\|
                |{Z}_s|^\delta  
            \big\|_{L^{16}(\Omega;\mathbb{R})}
            \big)
            \big(
                1+|x|
            \big)^{2c_g+1}
               h^{\tfrac{1}{2}},
    \end{aligned}
\end{equation}
which confirms the second assertion of \eqref{eq:ex_inte_property_condi_on_b}.
We mention that the assertions \eqref{eq:ex_inte_property_condi_on_a} and \eqref{eq:ex_inte_property_condi_on_b} for 
the stopped increment-tamed order $1.5$ method can be obtained in a similar manner as the stopped increment-tamed Milstein method without any difficulty. Thus we omit it here and the proof is completed.
\qed 

\section{Optimal strong convergence rates}
In this section we aim to recover optimal strong convergence rate of our new methods in a non-globally monotone setting.
%In light of the existing results \cite{Hutzenthaler2020,hutzenthaler2022stopped}, 
To be more concise, we restrictive ourselves to  convergence rates of stopped increment-tamed Milstein and order $1.5$ methods. 
The pathwise uniform strong convergence rate is considered here, which requires the diffusion term $g(x)$ to be Lipschitz. 
For the situation of non-Lipschitz diffusion term, the pointwise strong convergence rate can be obtained in a manner similar to \cite[Theorem 3.2, 1]{dai2023order}. 

For the simplicity of notations, we use the notations $\tilde{g}: \mathbb{R}^d \rightarrow \mathbb{R}$ and $\tilde{b}: \mathbb{R}^d \rightarrow \mathbb{R}$ to denote $g^{(ij)}$ and $b^{(ij)}$ for some fixed $i=1,...,d,j=1,...,m$.

\subsection{Convergence rate of stopped increment-tamed Milstein  method}

\begin{theorem}\label{thm:conver_rate_stopped_tamed_milstein_method}
     Let 
        $f\in 
         \mathcal{C}^1_{\mathcal{P}}(\mathbb{R}^d,\mathbb{R}^d) \cap {C}^2(\mathbb{R}^d,\mathbb{R}^d)$, let
         $g \in {C}^2(\mathbb{R}^d,\mathbb{R}^{d \times m})$ be Lipschitz 
    and let $g^{(j)\prime} \in \mathcal{C}^1_{\mathcal{P}}(\mathbb{R}^d,\mathbb{R}^{d\times d}),j=1,...,m$.
    Let $U_0\in \mathcal{C}^1_{\mathcal{P}}(\mathbb{R}^d,[0,\infty)) \cap {C}^2(\mathbb{R}^d,[0,\infty))$ with constants $K_{0},c_{0}\geq 0$ and 
    $U_1 \in C(\mathbb{R}^d,[0,\infty))$ satisfy $|U_1(x)|\leq K_{1}(1+|x|)^{c_{1}}, K_{1},c_{1}\geq 0$.
    Let the stopped increment-tamed Milstein method be defined by \eqref{eq:stop_tamed_milstein_method} with $\delta -2\theta \geq 3$ and 
        let $c,v,T\in(0,\infty),q,q_1,q_2\in (0,\infty],\alpha \in[0,\infty), p\geq 4 $.
        For any $x,y\in \mathbb R^d$, assume  that
        \begin{enumerate}[{\rm(1)}]
            \item 
                there exist constants $L,\kappa\geq0$ such that for any $i=1,...,d$,
                $$
                    \|\operatorname{Hess}_x(f^{(i)}(x))\|
                    \leq 
                    L(1+|x|)^{\kappa};
                $$
            \item
                $|x|^{1/c} \leq c(1+U_0(x))$\ and\ $\mathbb{E}\left[e^{U_0(X_0)}\right]<\infty;$
            \item
                $(\mathcal{A}_{f,g} U_0 )(x)+\tfrac{1}{2}|g(x)^{*}(\nabla U_0(x))|^2+U_1(x) \leq c+\alpha U_0(x);$
            \item
                $
                \langle x-y,f(x)-f(y)\rangle \leq \left[c + \tfrac{U_0(x)+U_0(y)}{2q_1Te^{\alpha T}}+\tfrac{U_1(x)+U_1(y)}{2q_2e^{\alpha T}}\right]|x-y|^2. 
                $
        \end{enumerate}
        Then for $\tfrac{1}{q}=\tfrac{1}{q_1}+\tfrac{1}{q_2},\tfrac{1}{v}=\tfrac{1}{p}+\tfrac{1}{q}$ and sufficiently small $h$, the scheme \eqref{eq:stop_tamed_milstein_method} 
        applied to SDEs \eqref{eq:typical_sde} is pathwise uniformly convergent with order one: 
        \begin{equation}
    \label{eq:thm_stopped_tamed_method_conclu001}
            \big\| \sup_{t\in[0,T]}|X_t-Y_t| \big\|_{L^{v}(\Omega;\mathbb{R})}
            \leq 
            C h.
        \end{equation}  
If the condition $(4)$ in Theorem {\rm \ref{thm:conver_rate_stopped_tamed_milstein_method}} is replaced by 
the following one: 
\begin{enumerate}[{\rm(4')}]
            \item 
for any $\eta>0$, there exists a constant $K_{\eta}$ such that  \begin{equation}\label{eq:strengthed_condition_item_stopped_tame}
            \langle x-y,f(x)-f(y)\rangle \leq \big[K_{\eta}+\eta\big(U_0(x)+U_0(y)+U_1(x)+U_1(y)\big)\big]|x-y|^2,
        \end{equation}
\end{enumerate}
        then  for any $v>0$ and sufficiently small $h$, it holds that
\begin{equation}\label{eq:strengthed_conclusion_stopped_tame}
            \big\| 
            \sup_{t\in[0,T]}|X_t-Y_t| 
            \big\|_{L^{v}(\Omega;\mathbb{R})}\leq Ch
            .
        \end{equation}  
\end{theorem}
\textbf{Proof:} We follow a treatment similar to  \cite[Theorem 4.2]{dai2023order} to prove the present theorem. Note that the moments of any order of exact solution and numerical solution are bounded,
due to condition (2), condition (3) and exponential integrability properties \eqref{eq:ex_inter_property_01}.
Thanks to equation \eqref{eq:ex_inter_property_01},   \cite[Theorem 3.2]{dai2023order} and some existing arguments in \cite[Theorem 4.2]{dai2023order},  it is only necessary to show that
$$
    \int_0^T
    \big\|
        \mathbbm{1}_{s\leq \tau_{N}}
        |
            {f(Y_{\lfloor s \rfloor}) - a(s)}
        |
    \big\|^2_{{L^{p}(\Omega;\mathbb{R})}}{\rm d}s
    \leq
    Ch^2
$$
and 
$$
    \int_0^T
    \big\|
        \mathbbm{1}_{s\leq \tau_{N}}
        |
            {g(Y_{s}) - b(s)}
        |
    \big\|^2_{{L^{p}(\Omega;\mathbb{R})}}{\rm d}s
    \leq
    Ch^2,
$$
where $C$ a positive constant independent of $h$.
For the first assertion, by Lemma \ref{lem:function_in_schme_property_02} and $\delta -2\theta \geq 3$,
\begin{equation}
        \begin{aligned}
        &
        \int_0^T
        \big\|
            \mathbbm{1}_{s\leq \tau_{N}}
            |
            {f(Y_{\lfloor s \rfloor}) - a(s)}
            |
        \big\|^2_{{L^{p}(\Omega;\mathbb{R})}}{\rm d}s
    \\
    &=
        \int_0^T
        \Big\|
            \mathbbm{1}_{s <  \tau_{N}}
            \big|
            (\Upsilon_h '({Z}_s)-I)
            f(Y_{\lfloor s \rfloor})
            +  
            \tfrac{1}{2} \sum\limits_{j=1}^m \Upsilon_h ''({Z}_s)
            \big(
                G_2(Y_{\lfloor s \rfloor},W_s)e_j,
                G_2(Y_{\lfloor s \rfloor},W_s)e_j
            \big)
            \big|
        \Big\|^2_{{L^{p}(\Omega;\mathbb{R})}}{\rm d}s
    \\
    &\leq 
    \int_0^T
        2\big\|
            |
            (\Upsilon_h '({Z}_s)-I)
            f(Y_{\lfloor s \rfloor})
            |
        \big\|^2_{{L^{p}(\Omega;\mathbb{R})}}
            + 
        \Big\|
            \big| 
             \sum\limits_{j=1}^m \Upsilon_h ''({Z}_s)
            \big(
               G_2(Y_{\lfloor s \rfloor},W_s)e_j,
               G_2(Y_{\lfloor s \rfloor},W_s)e_j
            \big)
            \big| 
        \Big\|^2_{{L^{p}(\Omega;\mathbb{R})}}{\rm d}s
    \\
    &\leq Ch^2, 
        \end{aligned}
\end{equation}
where 
$$
G_2(Y_{\lfloor s \rfloor},W_s)
=
        g(Y_{\lfloor s \rfloor})
        +
        \sum_{j=1}^{m}  
        \mathscr{L}_g^j g(Y_{\lfloor s \rfloor}) 
        \big(W^{(j)}_s-W^{(j)}_{\lfloor s \rfloor}\big),
\
{Z}_s=
        \int_{\lfloor s \rfloor}^s 
        f(Y_{\lfloor s \rfloor}) {\rm d}s
    +  
        \int_{\lfloor s \rfloor}^s 
         G_2(Y_{\lfloor s \rfloor},W_s)
          {\rm d}W_s.
$$
For the second assertion, applying It\^{o}'s formula to expand $\tilde{g}(Y_s)$ at $s=\lfloor s \rfloor$ yields
\begin{equation}\label{eq:g_Y_s_ito_expand} 
        \begin{aligned}
           \tilde{g}(Y_s)&=\tilde{g}({Y_{\lfloor s \rfloor}})
                        +
                        \int_{\lfloor s \rfloor}^s
                        \mathbbm{1}_{r< \tau_N}
                                \tilde{g}'(Y_r)
                                \Upsilon_h '(Z_{r})
                                g(Y_{\lfloor s \rfloor}) 
                                {\rm{d}}W_r
            \\
                    &\quad \quad
                    +
                    \int_{\lfloor s \rfloor}^s
                    \mathbbm{1}_{r< \tau_N}
                        \tilde{g}'(Y_r)
                        \Upsilon_h '(Z_{r})
                        \Big(
                        \sum_{j_1=1}^{m}  
                        \mathscr{L}_g^{j_1} g(Y_{\lfloor s \rfloor}) 
                        \big(W^{(j_1)}_r-W^{(j_1)}_{\lfloor s \rfloor}\big)
                        \Big) 
                        {\rm{d}}W_r
            \\
                    &\quad \quad
                    +
                    \int_{\lfloor s \rfloor}^s
                        \big(
                             \tilde{g}'(Y_r)a(r)
                             +
                            \tfrac{1}{2} \operatorname{tr}
                            \big(b(r)^* \operatorname{Hess}_x(\tilde{g}(Y_r)) b(r)\big)
                        \big)
                        {\rm{d}}r
            \\
                &=:
                        \tilde{g}({Y_{\lfloor s \rfloor}})
                        +
                        \int_{\lfloor s \rfloor}^s
                        \mathbbm{1}_{r< \tau_N}
                                \tilde{g}'(Y_r) 
                                \Upsilon_h '(Z_{r})
                                g(Y_{\lfloor s \rfloor}) 
                                {\rm{d}}W_r
                    +
                        R_M({\lfloor s \rfloor},s).
        \end{aligned}
    \end{equation}
Then by Lemma \ref{lem:function_in_schme_property_02} and $g^{(ij)\prime} \in \mathcal{C}^1_{\mathcal{P}}(\mathbb{R}^d,\mathbb{R})$
we obtain 
\begin{equation} 
        \begin{aligned}
        &\int_0^T
            \big\|
                \mathbbm{1}_{s\leq \tau_{N}}
                |
                {\tilde g(Y_{s}) - \tilde b(s)}
                |
            \big\|^2_{{L^{p}(\Omega;\mathbb{R})}}{\rm d}s
        \\
            &=
            \int_0^T
            \Big\|
                \mathbbm{1}_{s < \tau_{N}}
                \big|
                \tilde g(Y_{s}) 
                -
                \big(\Upsilon_h '(Z_{s})
                g(Y_{\lfloor s \rfloor})\big)^{(ij)} 
                -
                \big(
                        \sum_{j_1=1}^{m} 
                        \Upsilon_h '(Z_{s}) 
                        \mathscr{L}_g^{j_1} g(Y_{\lfloor s \rfloor}) 
                        (W^{(j_1)}_s-W^{(j_1)}_{\lfloor s \rfloor})
                \big)^{(ij)} 
                \big|
            \Big\|^2_{{L^{p}(\Omega;\mathbb{R})}}{\rm d}s
        \\
            &\leq 
            C
            \int_0^T
            \Big(
            \big\|
                    \big(
                    (\Upsilon_h '(Z_{s})-I)
                    g(Y_{\lfloor s \rfloor})
                    \big)^{(ij)}
            \big\|^2_{{L^{p}(\Omega;\mathbb{R})}}
            +
            \Big\|
                |
                    \int_{\lfloor s \rfloor}^{\tau_N}
                        %\mathbbm{1}_{r< \tau_N}
                                \tilde{g}'(Y_r) 
                                (\Upsilon_h '(Z_{r})-I)
                                g(Y_{\lfloor s \rfloor}) 
                                {\rm{d}}W_r
                |
            \Big\|^2_{{L^{p}(\Omega;\mathbb{R})}}
            \Big)
            {\rm d}s
        \\
            &\quad 
            +
            C
            \int_0^T
            \Big\|
                |
                    \int_{\lfloor s \rfloor}^s
                        \mathbbm{1}_{r< \tau_N}
                                \big(
                                \tilde{g}'(Y_r) 
                                    -\tilde{g}'(Y_{\lfloor s \rfloor})
                                \big)
                                g(Y_{\lfloor s \rfloor}) 
                                {\rm{d}}W_r
                |
            \Big\|^2_{{L^{p}(\Omega;\mathbb{R})}}
            {\rm d}s
        \\
            &\quad 
            +
            C
            \int_0^T
            \Big\|
                \big|
                \sum_{j_1=1}^m
                \big(
                        \tilde{g}'(Y_{\lfloor s \rfloor})
                        g^{(j_1)}(Y_{\lfloor s \rfloor})
                        -
                        (
                        \Upsilon_h '(Z_{s}) 
                        \mathscr{L}_g^{j_1} g(Y_{\lfloor s \rfloor})
                        )^{(ij)} 
                \big)
                        \big(W^{(j_1)}_s-W^{(j_1)}_{\lfloor s \rfloor}\big)
                \big|
            \Big\|^2_{{L^{p}(\Omega;\mathbb{R})}}{\rm d}s
                    \\
            &\quad 
            +
            C
            \int_0^T
            \big\|
            |
            R_M({\lfloor s \rfloor},s)
            |
            \big\|^2_{{L^{p}(\Omega;\mathbb{R})}}
            {\rm d}s
        \\
            &\leq 
            Ch^2.
        \end{aligned}
\end{equation}
The proof is thus completed.\qed

\subsection{Convergence rate of stopped increment-tamed order 1.5 method}

In order to show the expected strong convergence rate of the stopped increment-tamed order $1.5$ method, we start with the following essential lemma quoted from \cite[Lemma 3.1]{dai2023order}.
\begin{lemma}\label{lem:pre_estimate_theory}
    Let 
    $f:\mathbb R^d\rightarrow \mathbb R^d,g:\mathbb R^d\rightarrow \mathbb R^{d \times m} $ 
    be measurable functions. 
    Let $a: [0,T] \times \Omega \rightarrow \mathbb R^d,b: [0,T] \times \Omega \rightarrow \mathbb R^{d\times m}$ 
    be predictable stochastic processes and let $\tau:\Omega\rightarrow [0,T]$
    be a stopping time. 
    Let $\{ X_s \}_{s\in[0,T]}$ and 
    $\{ Y_s \}_{s\in[0,T]}$  be defined by \eqref{eq:typical_sde} and \eqref{eq:approxi_sde} with continuous sample paths, respectively.
    Assume that $\int_0^T |a(s)|+\|b(s)\|^2+|f(X_s)|+\|g(X_s)\|^2+|f(Y_s)|+\|g(Y_s)\|^2 {\rm{d}}s< \infty \ \mathbb{P}$-a.s. and for $\varepsilon \in (0,\infty),p\geq 2$ with $\mathbb{P}$-a.s.
        \begin{equation}\label{eq:assum_ex}
            \int_{0}^{\tau}\left[\tfrac{\left\langle X_{s}-Y_{s}, f\left(X_{s}\right)-f\left(Y_{s}\right)\right\rangle+\frac{(1+\varepsilon)(p-1)}{2}\left\|g\left(X_{s}\right)-g\left(Y_{s}\right)\right\|^{2}}{\left|X_{s}-Y_{s}\right|^{2}}\right]^{+} {\rm d} s<\infty.
        \end{equation}
    Then it holds
    %\mbox{}\vspace*{-1.8\baselineskip} 
\begin{align*}
            &
            \Big\|\sup_{t\in[0,T]} \tfrac{|X_{t \wedge \tau}-Y_{t \wedge \tau}|}{{\exp (\int_0^{t \wedge \tau} {\frac{1}{2}{\eta_{r}}} {\rm{d}}r)}}\Big\|_{{L^{p}(\Omega;\mathbb{R})}}
            \\
            &\leq \Bigg[
            \|{\xi_X-\xi_Y}{\|^2_{L^{p}(\Omega;\mathbb{R}^d)}}
            + 
            \underbrace{
            \Big\|\sup_{t\in[0,T]}\int_0^{t \wedge \tau} {\tfrac{{2\langle {X_s} - {Y_s},\left( {g({X_s}) -b(s)} \right){\rm{d}}W_s\rangle }}{{\exp (\int_0^s {{\eta_{r}}} {\rm{d}}r)}}}\Big\|_{{L^{p/2}(\Omega;\mathbb{R})}} 
            }_{=:\mathbb{T}_A}\quad \quad \quad
            \\
            &\ \ \ +\underbrace{
            \Big\|\sup_{t\in[0,T]}\int_0^{t \wedge \tau} {\tfrac{{2\langle {X_s} - {Y_s},{f({Y_s}) - a(s)} \rangle }}{{\exp (\int_0^s {{\eta_{r}}} {\rm{d}}r)}}} {\rm{d}}s\Big\|_{{L^{p/2}(\Omega;\mathbb{R})}}
            }_{=:\mathbb{T}_B}
            \\
            &\ \ \ +\underbrace{
            \Big\|\sup_{t\in[0,T]}\int_0^{t \wedge \tau} {\tfrac{{(1+\frac{1}{\varepsilon})\|g({Y_s}) - b(s){\|^2}}}{{\exp (\int_0^s {{\eta_{r}}} {\rm{d}}r)}}} {\rm{d}}s\Big\|_{{L^{p/2}(\Omega;\mathbb{R})}}
            }_{=:\mathbb{T}_D}
            \Bigg]^{1/2},
            \stepcounter{equation}
            \tag{\theequation}
            \label{eq:pre_estimate_with_sup}
\end{align*}
    where 
    \begin{equation}\label{eq:pre_estimate_eta}
        \eta_{r}:= 2 \mathbbm{1}_{r\leq \tau} (\omega)
        \left[\tfrac{\langle {X_r} - {Y_r},f(X_r)-f(Y_r)\rangle+{\frac{1+\varepsilon}{2}\|g(X_r)-g(Y_r){\|^2}
        }}{|{X_r} - {Y_r}|^2}\right]^{+}.
    \end{equation}  
\end{lemma}
Equipped with Lemma \ref{lem:pre_estimate_theory}, we are well-prepared to reveal the optimal  convergence rate of stopped increment-tamed order $1.5$ method.

\begin{theorem}\label{thm:conver_rate_stopped_tamed_1.5_order_method}
     Let 
        $f\in  {C}^3(\mathbb{R}^d,\mathbb{R}^d)$ and
        $\operatorname{Hess}_x(f^{(i)}) \in \mathcal{C}^1_{\mathcal{P}}(\mathbb{R}^d,\mathbb{R}^{d\times d})$,
        i=1,...,d.
        Let
         $g \in {C}^3(\mathbb{R}^d,\mathbb{R}^{d \times m})$ be Lipschitz
   and
   $\operatorname{Hess}_x(g^{(ij)}) \in \mathcal{C}^1_{\mathcal{P}}(\mathbb{R}^d,\mathbb{R}^{d\times d}),i=1,...,d,j=1,...,m$.
    Let $U_0\in \mathcal{C}^1_{\mathcal{P}}(\mathbb{R}^d,[0,\infty)) \cap {C}^2(\mathbb{R}^d,[0,\infty))$ with constants $K_{0},c_{0}\geq 0$ and $U_1 \in C(\mathbb{R}^d,[0,\infty))$ satisfy $|U_1(x)|\leq K_{1}(1+|x|)^{c_{1}}, K_{1},c_{1}\geq 0$.
    Let the stopped increment-tamed order $1.5$ method be defined by \eqref{eq:stop_tamed_1.5_order_method} with $\delta -2\theta \geq 4$ and 
        let $c,v,T\in(0,\infty),q,q_1,q_2\in (0,\infty],\alpha \in[0,\infty), p\geq 4 $.
        For any $x,y\in \mathbb R^d$, assume
        \begin{enumerate}[{\rm(1)}]
            % \item 
            %     there exist constants $L,\kappa\geq0$ such that for any $i=1,...,d$,
            %     $$
            %         \|\operatorname{Hess}_x(f^{(i)}(x))\|
            %         \leq 
            %         L(1+|x|)^{\kappa};
            %     $$
            \item
                $|x|^{1/c} \leq c(1+U_0(x))$\ and\ $\mathbb{E}\left[e^{U_0(X_0)}\right]<\infty;$
            \item
                $(\mathcal{A}_{f,g} U_0 )(x)+\tfrac{1}{2}|g(x)^{*}(\nabla U_0(x))|^2+U_1(x) \leq c+\alpha U_0(x);$
            \item
                $
                \langle x-y,f(x)-f(y)\rangle \leq \left[c + \tfrac{U_0(x)+U_0(y)}{2q_1Te^{\alpha T}}+\tfrac{U_1(x)+U_1(y)}{2q_2e^{\alpha T}}\right]|x-y|^2. 
                $
        \end{enumerate}
        Then for $\tfrac{1}{q}=\tfrac{1}{q_1}+\tfrac{1}{q_2},\tfrac{1}{v}=\tfrac{1}{p}+\tfrac{1}{q}$ and sufficiently small $h$, the approximation \eqref{eq:stop_tamed_milstein_method} 
        for \eqref{eq:typical_sde} admits 
        \begin{equation}
            \big\| \sup_{t\in[0,T]}|X_t-Y_t| \big\|_{L^{v}(\Omega;\mathbb{R})}
            \leq 
            C h^{\frac{3}{2}}.
        \end{equation}  
If the condition $(3)$ 
%in Theorem {\rm \ref{thm:conver_rate_stopped_tamed_1.5_order_method}} 
is replaced by the following one: 
\begin{enumerate}[{\rm(3')}]
            \item 
for any $\eta>0$, there exists a constant $K_{\eta}$ such that  \begin{equation}
            \langle x-y,f(x)-f(y)\rangle \leq \big[K_{\eta}+\eta\big(U_0(x)+U_0(y)+U_1(x)+U_1(y)\big)\big]|x-y|^2,
        \end{equation}
\end{enumerate}
        then  for any $v>0$ and sufficiently small $h$, it holds that
    \begin{equation}\label{eq:converg_rate_of_1.5_order_method}
            \big\|
                \sup_{t\in[0,T]}|X_t-Y_t|
            \big\|_{L^{v}(\Omega;\mathbb{R})}\leq Ch^{\frac{3}{2}}.
    \end{equation} 
\end{theorem}
\textbf{Proof:}
We employ Lemma \ref{lem:pre_estimate_theory} to show the desired assertions.
For $\mathbb{T}_A$, utilizing the Burkholder-Davis-Gundy inequality (see, e.g., \cite[Lemma 2.7]{wang2013tamed}), the Lipschitz assumption, the H{\"o}lder inequality and some elementary inequalities one arrives at
\begin{equation}\label{eq:estimate_T_A_with_sup}
\begin{aligned}
        \mathbb{T}_A
        &\leq 
            \Big\|
                \sup_{t\in[0,T]}
                \int_0^{t \wedge \tau_N} 
                    {\tfrac{{2\langle {X_s} - {Y_s},\left( {g({X_s})- g({Y_s})} \right){\rm{d}}W_s\rangle }}{{\exp (\int_0^s {{\eta_{r}}} 
                \rm{d}}r)}}
            \Big\|_{{L^{\frac{p}{2}}(\Omega;\mathbb{R})}}
            + 
            \Big\|
                \sup_{t\in[0,T]}
                \int_0^{t \wedge \tau_N} 
                    {\tfrac{{2\langle {X_s} - {Y_s},\left( g({Y_s})- b(s) \right){\rm{d}}W_s\rangle }}{{\exp (\int_0^s {{\eta_{r}}} 
                \rm{d}}r)}}
            \Big\|_{{L^{\frac{p}{2}}(\Omega;\mathbb{R})}}
        \\
            &\leq C_{g,p,m}
                \Big( 
                    \int_0^T
                    \Big\| 
                        {\tfrac{{\mathbbm{1}_{s\leq \tau_N}|{X_s} - {Y_s}|^2 }}
                        {{\exp (\int_0^s {{\eta_r}} {\rm{d}}r)}}}
                    \Big\|^2_{{L^{\frac{p}{2}}(\Omega;\mathbb{R})}}{\rm d}s
                \Big)^{\frac{1}{2}}
            +
            C_{p}
                \Big(
                    \int_0^T \sum\limits_{j=1}^{m}
                    \Big\| 
                    \mathbbm{1}_{s\leq \tau_N}
                        {\tfrac{|{X_s} - {Y_s}||g^{(j)}({Y_s})- b^{(j)}(s)| }
                        {\exp (\int_0^s {{\eta_{r}}
                        \rm{d}}r)}}
                    \Big\|_{{L^{\frac{p}{2}}(\Omega;\mathbb{R})}}
                    {\rm{d}}s
                \Big)^{\frac{1}{2}}
            \\
            &\leq 
                \tfrac{1}{4}
                \sup_{t\in[0,T]}
                \Big\| 
                    {\tfrac{{|{X_{t\wedge\tau_N}} - {Y_{t\wedge\tau_N}}| }}{{\exp (\int_0^{t\wedge\tau_N} {\frac{1}{2}{\eta_r}} {\rm{d}}r)}}}
                \Big\|^2_{{L^{p}(\Omega;\mathbb{R})}}
                +
                C\int_0^T  
                    \Big\|
                        {\tfrac{{{|{X_{s\wedge\tau_N}} -{Y_{s\wedge\tau_N}}|}}}{{\exp (\int_0^{s\wedge\tau_N} 
                        {\frac{1}{2}{\eta_r}} {\rm{d}}r)}}}
                    \Big\|^2_{{L^{p}(\Omega;\mathbb{R})}}
                 {\rm d}s
            \\
            &\quad
            +
            C\int_0^T
            \Big\|
                \mathbbm{1}_{s\leq \tau_N}\big({g({Y_s})-b(s)}\big) 
            \Big\|^2_{{L^{p}(\Omega;\mathbb{R}^{d\times m})}}
            {\rm d}s.
        \end{aligned}
\end{equation}
Concerning $\mathbb{T}_D $, one can easily see
\begin{equation} \label{eq:estimate_T_D_with_sup}
        \begin{aligned}
            \mathbb{T}_D 
            &\leq C
            \int_0^T 
                \big\|
                    \mathbbm{1}_{s\leq \tau_N}\big(g({Y_s}) - b(s)\big)
                \big\|^2_{{L^{p}(\Omega;\mathbb{R}^{d\times m})}} {\rm{d}}s.
        \end{aligned}
    \end{equation} 
We next claim that there exists a constant $C$ independent of $h$ such that 
\begin{equation}\label{eq:g_Y_s-b_sCh3} 
        \begin{aligned}
            \int_0^T 
                \big\|
                    \mathbbm{1}_{s\leq \tau_N}\big(g({Y_s}) - b(s)\big)
                \big\|^2_{{L^{p}(\Omega;\mathbb{R}^{d\times m})}} {\rm{d}}s
                \leq
                Ch^{3}.
        \end{aligned}
    \end{equation}
Recall that $ b(s)=\mathbbm{1}_{s< \tau_N} \Upsilon_h '(Z_{s})G_3(s,Y_{t_k},W_s)$ and the
notation $\tilde{g}({Y_s})=g^{(ij)}(Y_s)$ for some fixed $i=1,...,d, j=1,...,m$.
Expanding $\tilde{g}({Y_s})$ at $\lfloor s \rfloor$ results in  
    \begin{equation} 
        \begin{aligned}
            \tilde{g}({Y_s})&=\tilde{g}({Y_{\lfloor s \rfloor}})
                        +
                        \int_{\lfloor s \rfloor}^s
                        \Big(
                             \tilde{g}'(Y_r)a(r)
                             +
                            \tfrac{1}{2} \operatorname{tr}
                            \big(b(r)^* \operatorname{Hess}_x(\tilde{g}(Y_r)) b(r)\big)
                        \Big)
                        {\rm{d}}r
            \\
                    &\quad \quad
                    +
                    \sum_{j_1=1}^{m}
                    \int_{\lfloor s \rfloor}^s
                     \tilde{g}'(Y_r) b^{(j_1)}(r) 
                        {\rm{d}}W^{(j_1)}_r.
        \end{aligned}
    \end{equation}
One can thus split  $\tilde{g}({Y_s})-\tilde{b}(s)$ into three additional items:
    \begin{equation} 
        \begin{aligned}
            &\tilde{g}({Y_s})-\tilde{b}(s)
            \\
            &=
            \big(
            {g}({Y_{\lfloor s \rfloor}})
            -\mathbbm{1}_{s<\tau_N}
            \Upsilon_h'(Z_s)
            {g}({Y_{\lfloor s \rfloor}})
            \big)^{(ij)}
            \\
            &\quad  
            +
            \int_{\lfloor s \rfloor}^s
                    \Big(
                        \tilde{g}'(Y_r)a(r)
                        +
                        \tfrac{1}{2} \operatorname{tr}
                        \big(b(r)^* \operatorname{Hess}_x(\tilde{g}(Y_r)) b(r)\big)
                        -
                        \mathbbm{1}_{s<\tau_N}
                        \big(
                            \Upsilon_h' (Z_s)
                            \mathcal{A}_{f,g} {g}(Y_{\lfloor s \rfloor})
                        \big)^{(ij)}
                    \Big)
                {\rm{d}}r
            \\
            &\quad 
                +
                \sum_{j_1=1}^{m}
                \bigg[
                \int_{\lfloor s \rfloor}^s
                 \tilde{g}'(Y_r) b^{(j_1)}(r)
                 {\rm{d}}W^{(j_1)}_r
                -
                \mathbbm{1}_{s<\tau_N}
                \big(
                    \Upsilon_h'(Z_s)
                    \mathscr{L}_g^{(j_1)} {g}(Y_{\lfloor s \rfloor})
                \big)^{(ij)}
                \big(W^{(j_1)}_s-W^{(j_1)}_{\lfloor s \rfloor}\big)
            \\
            &\quad \quad \quad \quad \quad
                -
                \sum\limits_{j_2=1}^{m} 
                 \mathbbm{1}_{s<\tau_N}
                 \big(
                    \Upsilon_h'(Z_s)
                    \mathscr{L}_g^{j_2}
                    \mathscr{L}_g^{j_1}
                    {g}(Y_{\lfloor s \rfloor})
                \big)^{(ij)} 
                \int_{\lfloor s \rfloor}^s  
                \big(W^{(j_2)}_r-W^{(j_2)}_{\lfloor s \rfloor}\big)
                {\rm d}W^{(j_1)}_r
                \bigg]
            \\
            &=: A_1+A_2+A_3.
        \end{aligned}
    \end{equation}
Since $\delta-2\theta \geq 4$ and thanks to Lemma \ref{lem:function_in_schme_property_02}, it is clear to check
\begin{equation} 
        \begin{aligned}
            \int_0^T 
                \big\|
                    \mathbbm{1}_{s\leq \tau_N}
                    A_1
                \big\|^2_{{L^{p}(\Omega;\mathbb{R})}} {\rm{d}}s
                \leq
                Ch^{3}.
        \end{aligned}
\end{equation}
Again, using $\delta-2\theta \geq 4$, some elementary inequalities and the definitions of $a(r),b(r)$ and noting $\operatorname{Hess}_x(\tilde{g}) \in \mathcal{C}^1_{\mathcal{P}}(\mathbb{R}^d,\mathbb{R}^{d\times d})$ one infers
\begin{align*}
            & \int_0^T 
                \big\|
                    \mathbbm{1}_{s\leq \tau_N}
                    A_2
                \big\|^2_{{L^{p}(\Omega;\mathbb{R})}} {\rm{d}}s
            \\
            &
                \leq 
                Ch\int_0^T 
                    \int_{\lfloor s \rfloor}^s
                     \Big\|
                        \tilde{g}'(Y_r)a(r)
                        +
                        \tfrac{1}{2} \langle
                        b(r), 
                        \operatorname{Hess}_x(\tilde{g}(Y_r)) b(r)
                        \rangle_{HS}
                        -
                        \mathbbm{1}_{s<\tau_N}
                        \big(
                            \Upsilon_h' (Z_s)
                            \mathcal{A}_{f,g} {g}(Y_{\lfloor s \rfloor})
                        \big)^{(ij)}
                    \Big\|^2_{{L^{p}(\Omega;\mathbb{R})}}
                {\rm{d}}r
                {\rm{d}}s
            \\
            &\leq 
            Ch\int_0^T 
                    \int_{\lfloor s \rfloor}^s
                     \Big\|
                        \tilde{g}'(Y_r)
                        \Upsilon_h' (Z_r)
                        f(Y_{\lfloor s \rfloor})
                        +
                        \tfrac{1}{2} 
                        \langle
                        \Upsilon_h' (Z_r)g(Y_{\lfloor s \rfloor}), 
                        \operatorname{Hess}_x(\tilde{g}(Y_r))
                        \Upsilon_h' (Z_r)g(Y_{\lfloor s \rfloor}) 
                        \rangle_{HS}
            \\
            &\quad \quad \quad \quad \quad \quad \quad
                     -
                        \big(
                            \Upsilon_h' (Z_s)
                            \mathcal{A}_{f,g} {g}(Y_{\lfloor s \rfloor})
                        \big)^{(ij)}
                    \Big\|^2_{{L^{p}(\Omega;\mathbb{R})}}
                {\rm{d}}r
                {\rm{d}}s
        \\
        &\quad 
        +
        Ch\int_0^T 
                    \int_{\lfloor s \rfloor}^s
                     \Big\|
                        \tilde{g}'(Y_r)
                        \big(
                         a(r)
                            -
                         \mathbbm{1}_{r<\tau_N}
                          \Upsilon_h' (Z_r)
                            f(Y_{\lfloor s \rfloor})
                         \big)
        \\
        &\quad \quad \quad \quad \quad \quad \quad
                        +
                        \tfrac{1}{2} 
                        \mathbbm{1}_{r<\tau_N}
                        \big\langle
                        b(r)+\Upsilon_h' (Z_r)g(Y_{\lfloor s \rfloor}), 
                        \operatorname{Hess}_x(\tilde{g}(Y_r))
                        \big(
                        b(r)-\Upsilon_h' (Z_r)g(Y_{\lfloor s \rfloor})
                        \big) 
                        \big\rangle_{HS}
                        \Big\|^2_{{L^{p}(\Omega;\mathbb{R})}}
                {\rm{d}}r
                {\rm{d}}s
        \\
        &\leq 
            Ch\int_0^T 
                    \int_{\lfloor s \rfloor}^s
                     \Big\|
                        \tilde{g}'(Y_r)
                        (\Upsilon_h' (Z_r) -I)
                        f(Y_{\lfloor s \rfloor})
        \\
        &\quad \quad \quad \quad \quad \quad \quad
                        +
                        \tfrac{1}{2} 
                        \big \langle
                        \big(\Upsilon_h' (Z_r)+I\big)g(Y_{\lfloor s \rfloor}), 
                        \operatorname{Hess}_x(\tilde{g}(Y_r))
                        \big(\Upsilon_h' (Z_r)-I\big)g(Y_{\lfloor s \rfloor}) 
                        \big \rangle_{HS}
                        \Big\|^2_{{L^{p}(\Omega;\mathbb{R})}}
                {\rm{d}}r
                {\rm{d}}s
        \\
        &\quad 
        +
        Ch\int_0^T 
                    \int_{\lfloor s \rfloor}^s
                     \Big\|
                     \big(
                     \tilde{g}'(Y_r)
                     -
                     \tilde{g}'(Y_{\lfloor s \rfloor})
                     \big)
                        f(Y_{\lfloor s \rfloor})
                    +
                        \tfrac{1}{2} 
                        \big\langle
                        g(Y_{\lfloor s \rfloor}), 
                        \operatorname{Hess}_x
                        \big(   \tilde{g}(Y_r)-
                            \tilde{g}(Y_{\lfloor s \rfloor})
                        \big)
                        g(Y_{\lfloor s \rfloor}) 
                        \big\rangle_{HS}
                        \Big\|^2_{{L^{p}(\Omega;\mathbb{R})}}
                    {\rm{d}}r
                    {\rm{d}}s
        \\
        & \quad 
        +
        Ch\int_0^T 
                    \int_{\lfloor s \rfloor}^s
                     \Big\|
                     \tilde{g}'(Y_{\lfloor s \rfloor})
                        f(Y_{\lfloor s \rfloor})
                    +
                        \tfrac{1}{2} 
                        \langle
                        g(Y_{\lfloor s \rfloor}), 
                        \operatorname{Hess}_x
                        ( 
                            \tilde{g}(Y_{\lfloor s \rfloor})
                        )
                        g(Y_{\lfloor s \rfloor}) 
                        \rangle_{HS}
        \\
        &\quad \quad \quad \quad \quad \quad \quad
                    -
                        \big(
                            \Upsilon_h' (Z_s)
                            \mathcal{A}_{f,g} {g}(Y_{\lfloor s \rfloor})
                        \big)^{(ij)}
                    \Big\|^2_{{L^{p}(\Omega;\mathbb{R})}}
                {\rm{d}}r
                {\rm{d}}s 
            +Ch^3
        \\
        &
        \leq Ch^3.
         \stepcounter{equation}\tag{\theequation}
    \end{align*} 
For the estimate of $A_3$, according to  the definition of $b(r)$, some elementary inequalities, the application of It\^{o}'s formula to $\tilde{g}'(Y_r)$ at $r={\lfloor s \rfloor}$ and $\operatorname{Hess}_x(\tilde{g}) \in \mathcal{C}^1_{\mathcal{P}}(\mathbb{R}^d,\mathbb{R}^{d\times d})$,  we deduce that
\begin{align*}
            &\int_0^T 
                \big\|
                    \mathbbm{1}_{s \leq \tau_N}
                    A_3
                \big\|^2_{{L^{p}(\Omega;\mathbb{R})}} {\rm{d}}s
        \\ 
        &\leq
        C\sum^{m}_{j_1=1} 
                \int_0^T 
                    \Big\|
                    \mathbbm{1}_{s< \tau_N}
                    \Big[
                        \int_{\lfloor s \rfloor}^{s}
                            \tilde{g}'(Y_r)\Upsilon_h'(Z_r)
                            \big(
                                {g}^{(j_1)}(Y_{\lfloor s \rfloor})
                                +
                                \sum\limits_{j_2=1}^{m}  
                                \mathscr{L}_g^{j_2} g^{(j_1)}(Y_{\lfloor s \rfloor}) 
                                (W^{(j_2)}_r-W^{(j_2)}_{\lfloor s \rfloor})
                            \big) 
                            {\rm{d}}W^{(j_1)}_r
        \\
            &\quad \quad \quad \quad \quad
                -
                            \big(
                            \Upsilon_h'(Z_s)
                            \mathscr{L}_g^{(j_1)} {g}(Y_{\lfloor s \rfloor})
                            \big)^{(ij)}
                            \big(W^{(j_1)}_s-W^{(j_1)}_{\lfloor s \rfloor}\big)
         \\
            &\quad \quad \quad \quad \quad
                -
                \sum\limits_{j_2=1}^{m} 
                 \big(
                    \Upsilon_h'(Z_s)
                    \mathscr{L}_g^{j_2}
                    \mathscr{L}_g^{j_1}
                    {g}(Y_{\lfloor s \rfloor})
                \big)^{(ij)} 
                \int_{\lfloor s \rfloor}^s  
                \big(W^{(j_2)}_r-W^{(j_2)}_{\lfloor s \rfloor}\big)
                {\rm d}W^{(j_1)}_r
                \Big]
                \Big\|^2_{{L^{p}(\Omega;\mathbb{R})}}
                {\rm{d}}s
        \\ 
        &\quad 
        +
        C
                \int_0^T 
                    \Big\|
                        \int_{\lfloor s \rfloor}^{s}
                            \tilde{g}'(Y_r)\Upsilon_h'(Z_r)
                            \Big(
                                \mathcal{A}_{f,g} g(Y_{\lfloor s \rfloor}) (r-\lfloor s \rfloor)
        \\
        &\quad \quad \quad \quad \quad \quad 
                                +
                                \int_{\lfloor s \rfloor}^r
                                \sum\limits_{j_1,j_2=1}^{m}   
                                \mathscr{L}_g^{j_2}
                                \mathscr{L}_g^{j_1}
                                g(Y_{\lfloor s \rfloor}) 
                                (W^{(j_2)}_{r_1}-W^{(j_2)}_{\lfloor s \rfloor})
                                {\rm d}W^{(j_1)}_{r_1}
                                \Big)
                                {\rm{d}}W_r
                                \Big\|^2_{{L^{p}(\Omega;\mathbb{R}^d)}}
                                {\rm{d}}s
        \\
        &\leq 
        C\sum^{m}_{j_1=1} 
                \int_0^T 
                    \Big\|
                    \mathbbm{1}_{s< \tau_N}
                    \Big[
                        \int_{\lfloor s \rfloor}^{s}
                            \Big(
                            \tilde{g}'(Y_{\lfloor s \rfloor})
                                +
                                \int_{\lfloor s \rfloor}^r
                                \operatorname{Hess}_x(\tilde{g}(Y_{r_1}))
                                 b(r_1)
                                {\rm{d}}W_{r_1}
                            \Big)
                            \Upsilon_h'(Z_r)
        \\
            &\quad \quad \quad \quad \quad\quad \quad \quad \quad
            \cdot
                            \big(
                                {g}^{(j_1)}(Y_{\lfloor s \rfloor})
                                +
                                \sum\limits_{j_2=1}^{m}  
                                \mathscr{L}_g^{j_2} g^{(j_1)}(Y_{\lfloor s \rfloor}) 
                                (W^{(j_2)}_r-W^{(j_2)}_{\lfloor s \rfloor})
                            \big) 
                            {\rm{d}}W^{(j_1)}_r
        \\
            &\quad \quad \quad \quad \quad
                -
                            \big(
                            \Upsilon_h'(Z_s)
                            \mathscr{L}_g^{(j_1)} {g}(Y_{\lfloor s \rfloor})
                            \big)^{(ij)}
                            \big(W^{(j_1)}_s-W^{(j_1)}_{\lfloor s \rfloor}\big)
         \\
            &\quad \quad \quad \quad \quad
                -
                \sum\limits_{j_2=1}^{m} 
                 \big(
                    \Upsilon_h'(Z_s)
                    \mathscr{L}_g^{j_2}
                    \mathscr{L}_g^{j_1}
                    {g}(Y_{\lfloor s \rfloor})
                \big)^{(ij)} 
                \int_{\lfloor s \rfloor}^s  
                \big(W^{(j_2)}_r-W^{(j_2)}_{\lfloor s \rfloor}\big)
                {\rm d}W^{(j_1)}_r
                \Big]
                \Big\|^2_{{L^{p}(\Omega;\mathbb{R})}}
                {\rm{d}}s
                +Ch^3
        \\
        &\leq 
        C\sum^{m}_{j_1=1} 
                \int_0^T 
                    \Big\|
                        \int_{\lfloor s \rfloor}^{s}
                        \Big(
                            \tilde{g}'(Y_{\lfloor s \rfloor})
                             \Upsilon_h'(Z_r)
                                {g}^{(j_1)}(Y_{\lfloor s \rfloor})
                                +
                                \tilde{g}'(Y_{\lfloor s \rfloor})
                                \Upsilon_h'(Z_r)
                                \sum\limits_{j_2=1}^{m}  
                                \mathscr{L}_g^{j_2} g^{(j_1)}(Y_{\lfloor s \rfloor}) 
                                (W^{(j_2)}_r-W^{(j_2)}_{\lfloor s \rfloor})
        \\
            &\quad \quad \quad \quad \quad \quad\quad 
            +
                        \int_{\lfloor s \rfloor}^r
                                \operatorname{Hess}_x(\tilde{g}(Y_{r_1}))
                                 \Upsilon_h'(Z_{r_1})
                                 g(Y_{\lfloor s \rfloor})
                                {\rm{d}}W_{r_1}
                        \cdot
                        \Upsilon_h'(Z_r)
                        {g}^{(j_1)}(Y_{\lfloor s \rfloor})
                        \Big)
                        {\rm{d}}W^{(j_1)}_r
        \\
            &\quad \quad \quad \quad \quad
                -
                            \big(
                            \Upsilon_h'(Z_s)
                            \mathscr{L}_g^{(j_1)} {g}(Y_{\lfloor s \rfloor})
                            \big)^{(ij)}
                            \big(W^{(j_1)}_s-W^{(j_1)}_{\lfloor s \rfloor}\big)
         \\
            &\quad \quad \quad \quad \quad
                -
                \sum\limits_{j_2=1}^{m} 
                 \big(
                    \Upsilon_h'(Z_s)
                    \mathscr{L}_g^{j_2}
                    \mathscr{L}_g^{j_1}
                    {g}(Y_{\lfloor s \rfloor})
                \big)^{(ij)} 
                \int_{\lfloor s \rfloor}^s  
                \big(W^{(j_2)}_r-W^{(j_2)}_{\lfloor s \rfloor}\big)
                {\rm d}W^{(j_1)}_r
                \Big\|^2_{{L^{p}(\Omega;\mathbb{R})}}
                {\rm{d}}s
                +Ch^3
        \\
        &\leq 
        C\sum^{m}_{j_1=1} 
                \int_0^T 
                    \Big\|
                        \big(
                            \tilde{g}'(Y_{\lfloor s \rfloor})
                                {g}^{(j_1)}(Y_{\lfloor s \rfloor})
                                -
                            \big(
                            \Upsilon_h'(Z_s)
                            \mathscr{L}_g^{(j_1)} {g}(Y_{\lfloor s \rfloor})
                            \big)^{(ij)}
                        \big)
                        \big(
                        W_s^{(j_1)}-W_{\lfloor s \rfloor}^{(j_1)}
                        \big)
        \\
            &\quad \quad \quad \quad \quad 
            +
                \int_{\lfloor s \rfloor}^{s}
                \sum\limits_{j_2=1}^{m}
                \Big[ 
                \tilde{g}'(Y_{\lfloor s \rfloor})
                \mathscr{L}_g^{j_2} g^{(j_1)}(Y_{\lfloor s \rfloor}) 
            +
                \big(
                        \operatorname{Hess}_x(\tilde{g}(Y_{\lfloor s \rfloor}))
                        g(Y_{\lfloor s \rfloor})
                \big)^{(j_2)}
                {g}^{(j_1)}(Y_{\lfloor s \rfloor})
         \\
            &\quad \quad \quad \quad \quad \quad \quad
                - 
                 \big(
                    \Upsilon_h'(Z_s)
                    \mathscr{L}_g^{j_2}
                    \mathscr{L}_g^{j_1}
                    {g}(Y_{\lfloor s \rfloor})
                \big)^{(ij)} 
                 \Big]
                 \big(W^{(j_2)}_r-W^{(j_2)}_{\lfloor s \rfloor}\big)
                {\rm{d}}W^{(j_1)}_r
                \Big\|^2_{{L^{p}(\Omega;\mathbb{R})}}
                {\rm{d}}s
                +Ch^3
        \\
        &\leq Ch^3,
        \stepcounter{equation}\tag{\theequation}
\end{align*}
which confirms the assertion \eqref{eq:g_Y_s-b_sCh3}.
The estimate of $\mathbb{T}_B$ need to be treated more carefully.
First, expanding $f(Y_s)$ at $s={\lfloor s \rfloor}$ yields
\begin{equation}
    \begin{aligned}
        &f(Y_s)
    \\
        &=
        f(Y_{\lfloor s \rfloor})
        +
        \int_{\lfloor s \rfloor}^s 
        \big(
            f'(Y_r)a(r)+ 
            \tfrac{1}{2}
            \sum_{j=1}^m f''(Y_r)( b(r)e_j,b(r)e_j) 
        \big)
        {\rm d} r
        +
        \int_{\lfloor s \rfloor}^s 
            f'(Y_r)b(r)  
        {\rm d} W_r
    \\
        &=
        f(Y_{\lfloor s \rfloor})
        +
        \int_{\lfloor s \rfloor}^s 
        \mathbbm{1}_{r<\tau_N}
        \Big(
            f'(Y_{\lfloor s \rfloor})
            \Upsilon_h'(Z_r)
            f(Y_{\lfloor s \rfloor})
    \\
        &\quad \quad \quad \quad \quad \quad \quad
        + 
            \tfrac{1}{2}
            \sum_{j=1}^m f''(Y_{\lfloor s \rfloor})
            \big(  \Upsilon_h'(Z_r)g(Y_{\lfloor s \rfloor})e_j,  \Upsilon_h'(Z_r)g(Y_{\lfloor s \rfloor}) e_j
            \big) 
        \Big)
        {\rm d} r
    \\
        &\quad 
        +
         \int_{\lfloor s \rfloor}^s 
            \mathbbm{1}_{r<\tau_N} f'(Y_{\lfloor s \rfloor})
            \Upsilon_h'(Z_r)
            g(Y_{\lfloor s \rfloor}) 
            {\rm d} W_r
    \\
        &\quad 
        +
        \int_{\lfloor s \rfloor}^s 
        \mathbbm{1}_{r<\tau_N}
            \Big(
            \int_{\lfloor s \rfloor}^r 
                \operatorname{Hess}_x (f(Y_{r_1}))
                b(r_1)
                {\rm d} W_{r_1}
            \Upsilon_h'(Z_r)
            g(Y_{\lfloor s \rfloor})  
    \\
        &
        \quad \quad \quad \quad \quad \quad \quad
        +
        f'(Y_{\lfloor s \rfloor})
        \Upsilon_h'(Z_r)
        \sum_{j=1}^m 
        \mathscr{L}_g^j g(Y_{\lfloor s \rfloor})
            (W^{(j)}_r-W^{(j)}_{\lfloor s \rfloor})
        \Big)
        {\rm d} W_r
        +R_{1.5}({{\lfloor s \rfloor},s})
    \\
    & 
    =: B_1+ B_2+B_3+B_4+R_{1.5}({{\lfloor s \rfloor},s}),
    \end{aligned}
\end{equation}
where 
    \begin{align*}
        &R_{1.5}({{\lfloor s \rfloor},s})
    \\
        &:=
         \int_{\lfloor s \rfloor}^s 
        \Big(
        \big(
            f'(Y_{r})-f'(Y_{\lfloor s \rfloor})
        \big)a(r)
        + 
            \tfrac{1}{2}
            \sum_{j=1}^m 
            \big(
            f''(Y_r)-
            f''(Y_{\lfloor s \rfloor})
            \big)
            ( b(r)e_j, b(r)e_j ) 
        \Big)
        {\rm d} r
    \\
        &\quad 
        +
        \int_{\lfloor s \rfloor}^s 
        \Big(
        f'(Y_{\lfloor s \rfloor})
            (
                a(r)-
                \mathbbm{1}_{r<\tau_N}
                \Upsilon_h'(Z_r)
                f(Y_{\lfloor s \rfloor})
            )
    \\
        &\quad \quad \quad \quad   
        +  
            \tfrac{1}{2}
            \sum_{j=1}^m 
            f''(Y_{\lfloor s \rfloor})
            \big(
            \big( b(r)-
                \mathbbm{1}_{r<\tau_N}
                \Upsilon_h'(Z_r) g(Y_{\lfloor s \rfloor})\big)e_j,
                 \big(b(r)+\mathbbm{1}_{r<\tau_N}
                \Upsilon_h'(Z_r) g(Y_{\lfloor s \rfloor})\big) e_j 
            \big)
        \Big)
        {\rm d} r
    \\
        &\quad 
        +
        \int_{\lfloor s \rfloor}^s 
            \Big( 
                f'(Y_{r})
                -
                f'(Y_{\lfloor s \rfloor})
                -
                \int_{\lfloor s \rfloor}^r 
                \operatorname{Hess}_x (f(Y_{r_1}))
                b(r_1)
                {\rm d} W_{r_1}
            \Big)
            b(r)  
        {\rm d} W_r
    \\
        &\quad 
        +
         \int_{\lfloor s \rfloor}^s 
            \Big(
            \big(
            f'(Y_{\lfloor s \rfloor})
            +\int_{\lfloor s \rfloor}^r 
                \operatorname{Hess}_x (f(Y_{r_1}))
                b(r_1)
                {\rm d} W_{r_1}
            \big) 
    \\
    &\quad \quad \quad \quad  \quad 
            \cdot \big(
            b(r)
            - 
            \mathbbm{1}_{r<\tau_N}\Upsilon_h'(Z_r)
            g(Y_{\lfloor s \rfloor})
            -
            \mathbbm{1}_{r<\tau_N}\Upsilon_h'(Z_r)
            \sum\limits_{j=1}^{m}  
            \mathscr{L}_g^j g(Y_{\lfloor s \rfloor}) 
            (W^{(j)}_s-W^{(j)}_{\lfloor s \rfloor})
            \big)  
    \\
    &\quad \quad \quad \quad \quad  
            +
            \int_{\lfloor s \rfloor}^r 
                \operatorname{Hess}_x (f(Y_{r_1}))
                b(r_1)
                {\rm d} W_{r_1} 
                \cdot
            \mathbbm{1}_{r<\tau_N}
            \Upsilon_h'(Z_r)
                \sum\limits_{j=1}^{m}  
                \mathscr{L}_g^j g(Y_{\lfloor s \rfloor}) 
                (W^{(j)}_s-W^{(j)}_{\lfloor s \rfloor})  
        \Big) 
        {\rm d} W_r.
    \stepcounter{equation}\tag{\theequation}
    \end{align*}
Using It\^{o}'s formula to $f'(Y_r)$ and noting $\delta-2\theta \geq 4$, $\operatorname{Hess}_x(f^{(i)}) \in \mathcal{C}^1_{\mathcal{P}}(\mathbb{R}^d,\mathbb{R}^{d\times d})$ and Lemma \ref{lem:function_in_schme_property_02},  it is not difficult to check 
\begin{align*}
    &
    \Big\|
        \sup_{t\in[0,T]}
        \int_0^{t \wedge \tau_N}
        {\tfrac{{2\langle {X_s} - {Y_s},
        R_{1.5}({{\lfloor s \rfloor},s}) \rangle }}{{\exp (\int_0^s {{\eta_{r}}} {\rm{d}}r)}}} {\rm{d}}s
    \Big\|_{{L^{\frac{p}{2}}(\Omega;\mathbb{R})}}
    \leq 
    C
        \int_0^{T}
        \Big\|
        {\tfrac{{| {X_{s\wedge \tau_N}} - {Y_{s\wedge \tau_N}}|
         }}{{\exp (\int_0^{s\wedge \tau_N} {\frac12{\eta_{r}}} {\rm{d}}r)}}} 
         \Big\|^2_{{L^{p}(\Omega;\mathbb{R})}}
         {\rm{d}}s
    +Ch^3
    \stepcounter{equation}
    \tag{\theequation}
\end{align*}
and
\begin{equation} 
        \begin{aligned}
        &
            \Big\|
        \sup_{t\in[0,T]}
        \int_0^{t \wedge \tau_N}
        {\tfrac{{2\langle {X_s} - {Y_s},
       B_1+B_2+B_3 -a(s)
        \rangle }}{{\exp (\int_0^s {{\eta_{r}}} {\rm{d}}r)}}} {\rm{d}}s
            \Big\|_{{L^{\frac{p}{2}}(\Omega;\mathbb{R})}} \leq 
             C
        \int_0^{T}
        \Big\|
        {\tfrac{{| {X_{s\wedge \tau_N}} - {Y_{s\wedge \tau_N}}|
         }}{{\exp (\int_0^{s\wedge \tau_N} {\frac12{\eta_{r}}} {\rm{d}}r)}}} 
         \Big\|^2_{{L^{p}(\Omega;\mathbb{R})}}
         {\rm{d}}s
         +Ch^3.
        \end{aligned}
\end{equation}
Furthermore, by $\delta-2\theta \geq 4$, Lemma \ref{lem:function_in_schme_property_02} and H{\"o}lder's inequality, it holds 
\begin{equation} 
        \begin{aligned}
        &
        \Big\|
        \sup_{t\in[0,T]}
        \int_0^{t \wedge \tau_N}
        {\tfrac{{2\langle {X_s} - {Y_s},
       B_4
        \rangle }}{{\exp (\int_0^s {{\eta_{r}}} {\rm{d}}r)}}} {\rm{d}}s
            \Big\|_{{L^{p/2}(\Omega;\mathbb{R})}}
        \\
        &\leq 
        \Big\|
        \sup_{t\in[0,T]}
        \int_0^{t \wedge \tau_N}
        {\tfrac{{2 
        \left \langle 
        {X_s} - {Y_s},
       \int_{\lfloor s \rfloor}^s 
            \sum_{j=1}^m
            \int_{\lfloor s \rfloor}^r 
            \left(
                (\operatorname{Hess}_x (f(Y_{r_1}))
                b(r_1))^{(j)}
                g(Y_{\lfloor s \rfloor}) 
            +
            f'(Y_{\lfloor s \rfloor}) 
            \mathscr{L}_g^j g(Y_{\lfloor s \rfloor})
            \right)
            {\rm d} W^{(j)}_{r_1}
            {\rm d} W_r
        \right \rangle }}{{\exp (\int_0^s {{\eta_{r}}} {\rm{d}}r)}}} {\rm{d}}s
            \Big\|_{{L^{\frac{p}{2}}(\Omega;\mathbb{R})}}
        \\
        &\quad
        +
        C
        \int_0^{T}
        \Big\|
        {\tfrac{{| {X_{s\wedge \tau_N}} - {Y_{s\wedge \tau_N}}|
         }}{{\exp (\int_0^{s\wedge \tau_N} {\frac12{\eta_{r}}} {\rm{d}}r)}}} 
         \Big\|^2_{{L^{p}(\Omega;\mathbb{R})}}
         {\rm{d}}s
    +Ch^3.
        \end{aligned}
\end{equation}
As a consequence, one derives
\begin{align*}
    \mathbb{T}_B   
        &\leq 
        \Big\|
        \sup_{t\in[0,T]}
        \int_0^{t \wedge \tau_N}
        {\tfrac{{2 
        \left \langle 
        {X_s} - {Y_s},
       \int_{\lfloor s \rfloor}^s 
            \sum_{j=1}^m
            \int_{\lfloor s \rfloor}^r 
            \left(
                (\operatorname{Hess}_x (f(Y_{r_1}))
                b(r_1))^{(j)}
                g(Y_{\lfloor s \rfloor}) 
            +
            f'(Y_{\lfloor s \rfloor}) 
            \mathscr{L}_g^j g(Y_{\lfloor s \rfloor})
            \right)
            {\rm d} W^{(j)}_{r_1}
            {\rm d} W_r
        \right \rangle }}{{\exp (\int_0^s {{\eta_{r}}} {\rm{d}}r)}}} {\rm{d}}s
            \Big\|_{{L^{\frac{p}{2}}(\Omega;\mathbb{R})}}
        \\
        &\quad
        +
        C
        \int_0^{T}
        \Big\|
        {\tfrac{{| {X_{s\wedge \tau_N}} - {Y_{s\wedge \tau_N}}|
         }}{{\exp (\int_0^{s\wedge \tau_N} {\frac12{\eta_{r}}} {\rm{d}}r)}}} 
         \Big\|^2_{{L^{p}(\Omega;\mathbb{R})}}
         {\rm{d}}s
    +Ch^3.
        \stepcounter{equation}
        \label{eq:estimate_1.5_rate_T_B_temp}
        \tag{\theequation}
\end{align*}
The first term of equation \eqref{eq:estimate_1.5_rate_T_B_temp} can be estimated similarly to \cite[Theorem 3.2, (3.31)-(3.39)]{dai2023order}, due to the facts that
$$
\big\|
    |f(Y_r)-a(r)|
\big\|_{{L^{p}(\Omega;\mathbb{R})}}
\leq Ch,\ r \in [\lfloor s \rfloor,s]
$$
and assertion \eqref{eq:g_Y_s-b_sCh3}, which are used in the estimation of $\widetilde{B}_2$ and $\widetilde{B}_3$ in \cite[Theorem 3.2]{dai2023order}.  We hence arrive at
\begin{align*}
    \mathbb{T}_B   
        &\leq 
        \tfrac{1}{4}
        \sup_{s\in[0,T]}
        \Big\|
        {\tfrac{{| {X_{s\wedge \tau_N}} - {Y_{s\wedge \tau_N}}|
         }}{{\exp (\int_0^{s\wedge \tau_N} {\frac12{\eta_{r}}} {\rm{d}}r)}}} 
         \Big\|^2_{{L^{p}(\Omega;\mathbb{R})}}
        +
        C
        \int_0^{T}
        \Big\|
        {\tfrac{{| {X_{s\wedge \tau_N}} - {Y_{s\wedge \tau_N}}|
         }}{{\exp (\int_0^{s\wedge \tau_N} {\frac12{\eta_{r}}} {\rm{d}}r)}}} 
         \Big\|^2_{{L^{p}(\Omega;\mathbb{R})}}
         {\rm{d}}s
    +Ch^3.
        \stepcounter{equation}
        \tag{\theequation}
\end{align*}
Applying  the Gronwall inequality one obtains that
\begin{equation}
    \Big\|
        \sup_{t\in[0,T]} \tfrac{|X_{t \wedge \tau_N}-Y_{t \wedge \tau_N}|}{{\exp (\int_0^{t \wedge \tau_N} {\frac{1}{2}{\eta_{r}}} {\rm{d}}r)}}
    \Big\|^2_{{L^{p}(\Omega;\mathbb{R})}}
    \leq Ch^3,
\end{equation}
which gives, for $\tfrac{1}{p}+\tfrac{1}{q}=\tfrac{1}{v}$, 
\begin{equation}
        \begin{aligned}
            &\big\|\sup_{t\in[0,T]} |X_{t \wedge \tau_N}-Y_{t \wedge \tau_N}|\big\|_{{L^{v}(\Omega;\mathbb{R})}}
            \\
            &\leq    \Big\|\sup_{t\in[0,T]}\tfrac{|X_{t \wedge \tau_N} - Y_{t \wedge \tau_N}|}{\exp (\int_0^{t \wedge \tau_N} {\frac{1}{2}{\eta_r}} {\rm d}r)} \Big\|_{L^{p}(\Omega;\mathbb{R})}\cdot \Big\| {\exp \Big(\int_0^{\tau_N} {\tfrac{1}{2}{\eta_r}} {\rm{d}}r\Big)}\Big\|_{{L^{q}(\Omega;\mathbb{R})}}
            \\
            &\leq Ch^{\frac{3}{2}}.
        \end{aligned}
\end{equation}
Using exponential integrability property
\eqref{eq:ex_inter_property_01} and 
following a similar argument as in \cite[Theorem 4.2, (4.26)-(4.29)]{dai2023order}, we finally finish the proof.
\qed

\section{Examples and numerical experiments}

In this section we list several practical models 
that can fall within the above theoretical framework,
by choosing proper control functions $U_0(x)$ and $U_1(x)$. Regarding the choices of $U_0$ and $U_1$ for specific models, one can refer to \cite[3.12--3.17]{Hutzenthaler2020} and \cite[Theorem 5.6] {dai2023order} for more details.
Then we apply the above theoretical results to these models and give some numerical experiments to support the theoretical findings. For simplicity, the initial value $X_0$ of the following models are always assumed to be deterministic. 
%Moreover, for the choices of the control functions $U_0(x)$ and $U_1(x)$ of specific models, one can refer to \cite[3.12--3.17]{Hutzenthaler2020} and \cite[Theorem 5.6] {dai2023order}. 

\textbf{Stochastic Lorenz equation with additive noise.}\ Let $d=m=3$ and $\alpha_1,\alpha_2,\alpha_3 \geq 0$. For $x=(x_1,x_2,x_3)^*\in \mathbb{R}^3$, let
    \begin{equation}
        f(x)=(\alpha_1(x_2-x_1),\alpha_2x_1-x_2-x_1x_3,x_1x_2-\alpha_3x_3)
    \end{equation} 
    and $g(x)$ be some constant matrix.

\textbf{Brownian dynamics.}\ Let $d=m\geq1,c,\beta > 0$ and $\lambda \in [0,\tfrac{2}{\beta})$. 
    Assume that $V\in  C^3\big(\mathbb{R}^d,[0,\infty)\big)$, $\operatorname{Hess}_x \big(V^{(i)}\big)\in \mathcal{C}^1_{\mathcal{P}}\big(\mathbb{R}^d,\mathbb{R}^{d\times d}\big),i=1,...,d$ and 
    $\limsup_{r\searrow0} \sup_{z\in \mathbb{R}^d}
    $ $\tfrac{|z|^r}{1+V(z)}<\infty.$ For $x=(x_1,...,x_d)^*\in \mathbb{R}^d$, let
    \begin{equation}
        f(x)=-(\nabla V)(x),
        \ 
        g(x)=\sqrt{\beta}I_{\mathbb{R}^{d\times d}}.
    \end{equation} 
    Besides, we suppose that $(\Delta V)(x)\leq c+cV(x)+\lambda\|(\nabla V)(x)\|^2$ and for any $\eta>0$
    \[
        \sup_{x,y\in \mathbb{R}^d,x\neq y} \Big[\tfrac{\langle x-y,(\nabla V)(y)-(\nabla V)(x)\rangle}{|x-y|^2}-\eta\big(V(x)+V(y)+|(\nabla V)(x)|^2+|(\nabla V)(y)|^2\big)\Big]<\infty.
    \]

\textbf{Langevin dynamics.}\ Let $d=2m\geq1,\gamma\geq 0$ and $\beta>0$. 
    Assume that $V\in C^3\big(\mathbb{R}^m,[0,\infty)\big)$, $\operatorname{Hess}_x \big(V^{(i)}\big)\in \mathcal{C}^1_{\mathcal{P}}\big(\mathbb{R}^m,\mathbb{R}^{m\times m}\big),i=1,...,m$ and 
    $\limsup$ $_{r\searrow 0}
    $ $\sup_{z\in \mathbb{R}^m}\tfrac{|z|^r}{1+V(z)}<\infty.$ For $x=(x_1,x_2)^*\in \mathbb{R}^{2m},u\in\mathbb{R}^{m}$, let
    \begin{equation} \label{eq:Langevin_dynamics}
        f(x)=(x_2,-(\nabla V)(x_1)-\gamma x_2),g(x)u=(0,\sqrt \beta u).
    \end{equation} 
    In addition, we suppose that for any $\eta>0$
    \[
        \sup_{x,y\in \mathbb{R}^m,x\neq y} \Big[\tfrac{|(\nabla V)(x)-(\nabla V)(y)|}{|x-y|}-\eta\big(V(x)+V(y)+|x|^2+|y|^2\big)\Big]<\infty.
    \] 

\textbf{Experimental psychology model.}
Let $d=2,m=1$ and $\gamma,\lambda >0, \beta \in \mathbb{R}$. For $x=(x_1,x_2)^* \in \mathbb{R}^2$, let
$$
f(x)=\big(x_2^2(\lambda+4\gamma x_1)-\tfrac{1}{2}\beta^2 x_1,
,
-x_1x_2(\lambda+4\gamma x_1)-\tfrac{1}{2}\beta^2 x_2
\big)
, \
g(x)=(-\beta x_2, \beta x_1).
$$

    \textbf{Stochastic van der Pol oscillator.}\ Let $d=2,m\geq1, c,\lambda >0$ and $\gamma,\beta \geq 0$. For $x=(x_1,x_2)^*\in \mathbb{R}^2,u \in \mathbb{R}^m$, define
    \begin{equation}
        f(x)=(x_2,(\gamma-\lambda x^2_1)x_2-\beta x_1)^*,\ g(x)u=(0,\phi(x_1)u)^*,
    \end{equation}  
    where $\phi\in C^{2}(\mathbb{R}; \mathbb{R}^{m})$ is a globally Lipschitz function.

\textbf{Stochastic Duffing-van der Pol oscillator.}\ Let $d=2,m\geq1, \alpha_1,\alpha_2\in \mathbb{R}$ and $\alpha_3,c>0$. For $x=(x_1,x_2)^*\in \mathbb{R}^2,u \in \mathbb{R}^m$, define
    \begin{equation}
        f(x)=(x_2,\alpha_2x_2-\alpha_1x_1-\alpha_3x^2_1x_2-x^3_1)^*,\ g(x)u=(0,\phi(x_1)u)^*,
    \end{equation} 
    where $\phi\in C^{2}(\mathbb{R}; \mathbb{R}^{m})$ is a globally Lipschitz function. 

\textbf{Stochastic Lotka-Volterra competition model.}\ 
Let $b\in \mathbb{R}^d,\sigma \in \mathbb{R}^{d\times m}$ and $A=(a^{(ij)}) \in \mathbb{R}^{d\times d}, i,j=1,...,d$. Assume that every element of $A$ is non-negative and $\min_{1\leq i\leq d}\{a^{(ii)}\}>0$. For $x\in \mathbb{R}^d$, let
    \begin{equation}
        f(x)={\rm diag}(x)(b-Ax),\ g(x)u={\rm diag}(x)\sigma, \ X_0 \in \{x: x^{(1)}>0,...,x^{(d)}>0\},
    \end{equation} 
    where ${\rm diag}(x)$ represents a $d\times d$ diagonal matrix with principle diagonal $x$.

Applying the stopped increment-tamed Milstein method \eqref{eq:stop_tamed_milstein_method} with $\delta -2\theta \geq 3$ to all the models presented above, one obtains that 
$$
            \big\| \sup_{t\in[0,T]}|X_t-Y_t| \big\|_{L^{v}(\Omega;\mathbb{R})}
            \leq 
            C h
$$
for any $v>0$ and small $h$. Furthermore, if the stopped increment-tamed order $1.5$ method \eqref{eq:stop_tamed_1.5_order_method} is used to approximate  these models with $\delta -2\theta \geq 4$, then 
$$
            \big\| \sup_{t\in[0,T]}|X_t-Y_t| \big\|_{L^{v}(\Omega;\mathbb{R})}
            \leq 
            C h^{\frac{3}{2}}
$$
for any $v>0$ and small $h$. We mention that the potential function $V(x)$ in the Brownian dynamics and the Langevin dynamics and the diffusion coefficient function $\phi(x)$ in the stochastic van der Pol oscillator stochastic Duffing-van der Pol oscillator 
are required to be smoother, in order to obtain higher order like order $1.5$ (compare Theorem \ref{thm:conver_rate_stopped_tamed_milstein_method} with Theorem \ref{thm:conver_rate_stopped_tamed_1.5_order_method}). 

\begin{figure}[H]
      \centering
      \includegraphics[scale=0.2]{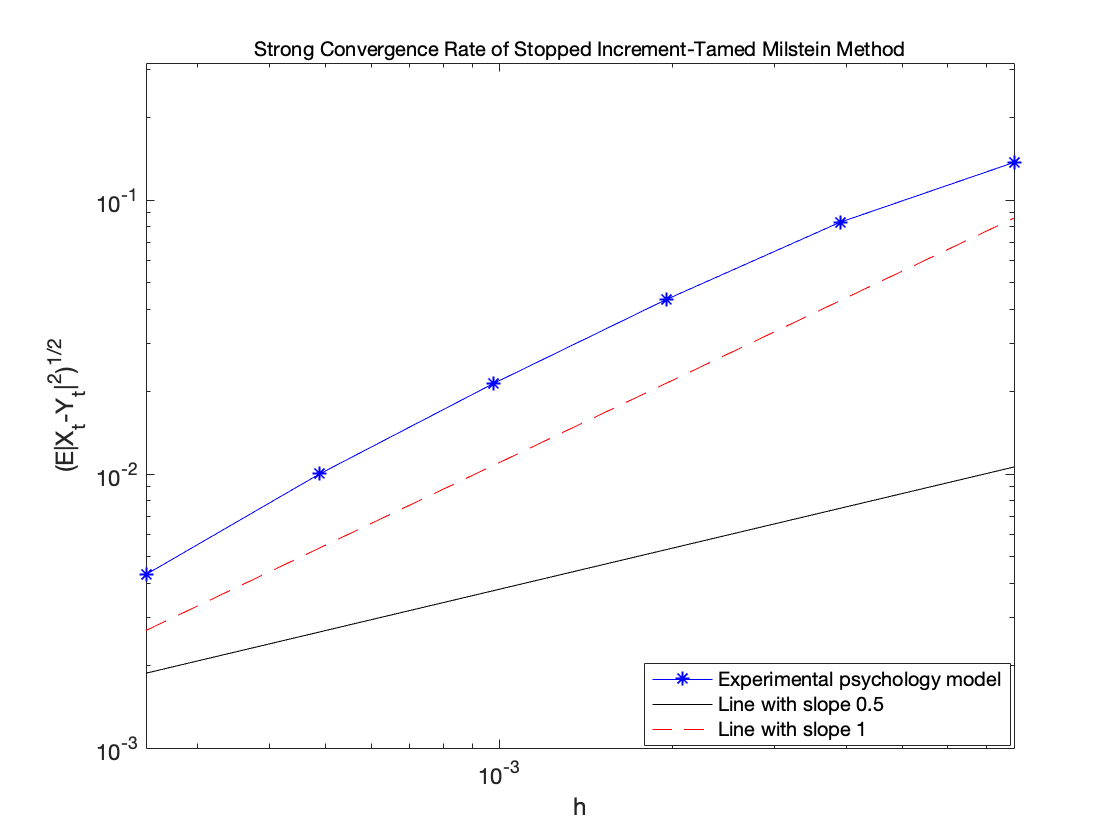}
      \caption{\small  Convergence rates
      of stopped increment-tamed Milstein method
      \label{fig:convergence_rate_01}}
    \end{figure}

\begin{figure}[H]
      \centering
      \includegraphics[scale=0.2]{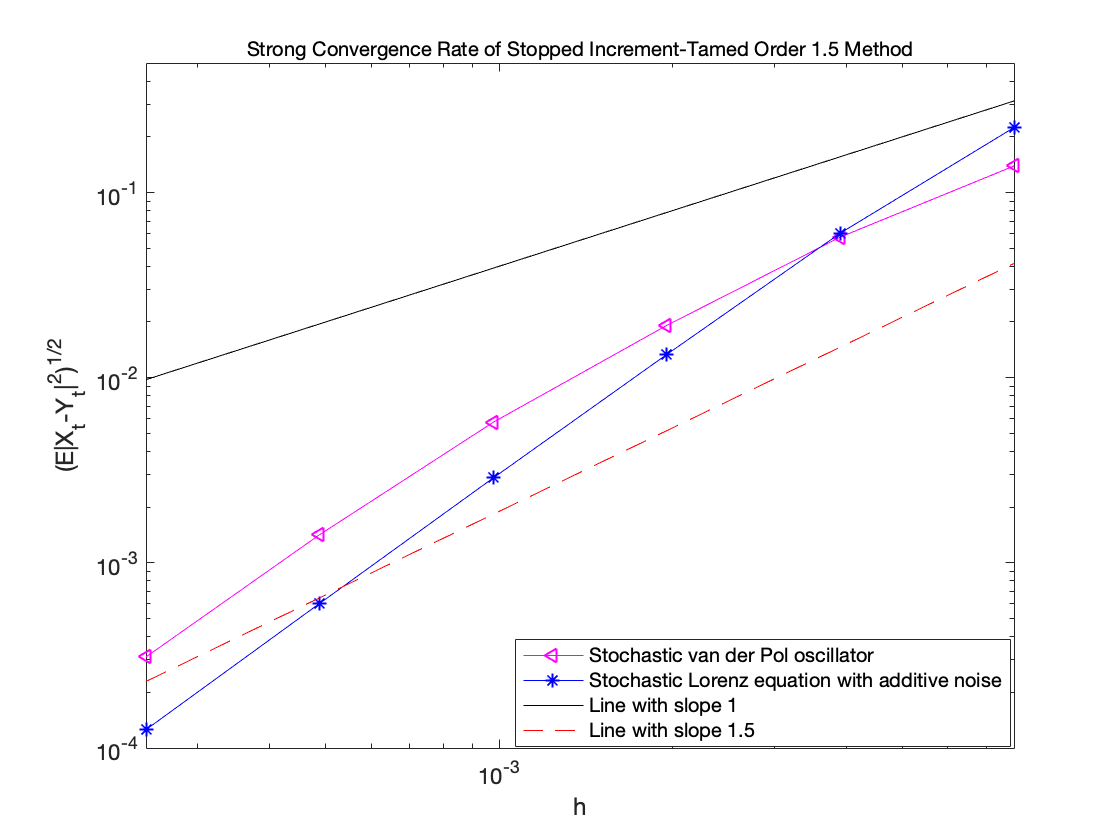}
      \caption{\small Convergence rates
      of stopped increment-tamed order $1.5$ method
      \label{fig:convergence_rate_02}}
    \end{figure}

%It is worthwhile to point out that our results significantly improve the existing works. 
%
Finally, we test the above theoretical results by performing some numerical results. To simplify the presentation, we take the stochastic Lorenz equation with additive noise, the experimental psychology model and the stochastic van der Pol oscillator for examples. 
Let $T = 1$, $N = 2^{k}, k = 7, 8,..., 12$ and regard fine approximations with $N_{\text{exact}} = 2^{14}$ as the "true" solution. The mean-square errors is considered and the average over $M = 5000$ Monte Carlo sample paths is used to approximate the expectation. By setting
\begin{equation}\label{eq:numer_test_method_para}
\delta=5, \theta=\tfrac{1}{4}, \gamma_1=1,  \gamma_2=1, \gamma_3=0.5,
\end{equation}
we plot the approximation errors of the stopped increment-tamed Milstein method
in Figure \ref{fig:convergence_rate_01} and those of the stopped increment-tamed order $1.5$ method
in Figure \ref{fig:convergence_rate_02},
where 
one can clearly detect the strong convergence rates of order $1$ and order $1.5$, respectively. This confirms the above theoretical findings.

\bibliographystyle{abbrv}

%\bibliography{../bib/bibfile}
\bibliography{Higher_order_reference}

\end{document}